\input amstex
\input amsppt.sty
\magnification=\magstep1
\hsize=33truecc
\vsize=22.2truecm
\baselineskip=16truept
\NoBlackBoxes
\TagsOnRight \pageno=1 \nologo
\def\Z{\Bbb Z}
\def\N{\Bbb N}

\def\Q{\Bbb Q}

\def\l{\left}
\def\r{\right}
\def\bg{\bigg}
\def\({\bg(}
\def\[{\bg\lfloor}
\def\){\bg)}
\def\]{\bg\rfloor}
\def\t{\text}
\def\f{\frac}

\def\sm{\setminus}

\def\bi{\binom}
\def\eq{\equiv}

\def\ls{\leqslant}
\def\gs{\geqslant}
\def\mo{\roman{mod}}

\def\Proof{\noindent{\it Proof}}

\def\Remark{\medskip\noindent{\it  Remark}}

\def\Ack{\medskip\noindent {\bf Acknowledgments}}
\hbox {{\it Number Theory and Related Area} (eds., Y. Ouyang, C.
Xing, F. Xu and P. Zhang),} \hbox{Adv. Lect. Math. 27, Higher
Education Press and International Press,} \hbox {Beijing--Boston,
2013, pp. 149--197.}
\bigskip
\topmatter
\title Conjectures and results on $x^2\ \mo\ p^2$ with $4p=x^2+dy^2$ \endtitle
\author Zhi-Wei Sun\endauthor
\date {\it In honor of Prof. Keqin Feng on the occasion of his 70th birthday}\enddate
\leftheadtext{Zhi-Wei Sun}
\rightheadtext{Conjectures and results on $x^2\ \mo\ p^2$ with $4p=x^2+dy^2$}
\affil Department of Mathematics, Nanjing University\\
 Nanjing 210093, People's Republic of China
  \\  zwsun\@nju.edu.cn
  \\ {\tt http://math.nju.edu.cn/$\sim$zwsun}
\endaffil
\abstract Given a squarefree positive integer $d$, we want to
find {\it integers} (or rational numbers with denominators not divisible by large primes)
 $a_0,a_1,a_2,\ldots$ such that for sufficiently large primes $p$ we have
  $\sum_{k=0}^{p-1}a_k\eq x^2-2p\pmod{p^2}$
 if $4p=x^2+dy^2$ (and $4\nmid x$ if $d=1$), and $\sum_{k=0}^{p-1}a_k\eq0\pmod{p^2}$ if $(\f{-d}p)=-1$.
We give a survey of conjectures and results on this topic and point out the connection between this problem
and series for $1/\pi$.
\endabstract
\thanks 2010 {\it Mathematics Subject Classification}.\,Primary 11E25, 11A07;
Secondary  05A10, 11B65.
\newline\indent {\it Keywords}. Representations of primes by binary quadratic forms, series for $1/\pi$.
\newline\indent Supported by the National Natural Science
Foundation (grant 11171140) of China.
\endthanks
\endtopmatter
\document

\heading{1. Introduction}\endheading

Let $d\in\Z^+=\{1,2,3,\ldots\}$ be squarefree.
If an odd prime $p$ not dividing $d$ can be
written in the form $x^2+dy^2$ with $x,y\in\Z$
then the Legendre symbol $(\f{-d}p)$ must be equal to 1. For $d=1,2,3$, it is well known that
any odd prime $p$ with $(\f{-d}p)=1$ can be written uniquely in the form $x^2+dy^2$ with $x,y\in\Z^+$
(and $2\nmid x$ if $d=1$). (See, e.g., Cox [C].)

Let $p\eq1\pmod4$ be a prime and write $p=x^2+y^2$ with $x\eq1\pmod4$ and $y\eq0\pmod2$.
In 1828 Gauss determined $x$ mod $p$ by showing the congruence
$$\bi{(p-1)/2}{(p-1)/4}\eq2x\pmod p.$$
In 1986, S. Chowla, B. Dwork and R. J. Evans [CDE] used Gauss and Jacobi sums to prove that
$$\bi{(p-1)/2}{(p-1)/4}\eq\f{2^{p-1}+1}2\l(2x-\f p{2x}\r)\pmod {p^2},$$
which was first conjectured by F. Beukers. This implies that
$$\bi{(p-1)/2}{(p-1)/4}^2\eq 2^{p-1}(4x^2-2p)\pmod{p^2}$$
since
$$\l(2x-\f p{2x}\r)^2\eq 4x^2-2p\pmod{p^2}.$$
In 2010 J. B. Cosgrave and K. Dilcher [CD] determined $\bi{(p-1)/2}{(p-1)/4}$ mod $p^3$.
A result due to Cauchy and Whiteman (cf. [HW84]) asserts that if $p\eq1\pmod{20}$ then
$$\bi{(p-1)/2}{(p-1)/20}\eq(-1)^{[5\mid x]}\bi{(p-1)/2}{3(p-1)/20}\pmod p.$$
(Throughout this paper, for an assertion $A$ we let $[A]$ takes $1$ or $0$ according as $A$ holds or not.)
The reader may consult [HW79] for more such congruences of the Cauchy type.

In 1989 K. M. Yeung [Ye] employed the Gross-Koblitz formula
([GK]) for the $p$-adic Gamma function to prove that if $p\eq1\pmod6$ is a prime and $p=x^2+3y^2$ with $x,y\in\Z$
and $3\mid x-1$ then
$$\bi{(p-1)/2}{(p-1)/6}\eq\l(2x-\f p{2x}\r)\l(1-\f{2^p-2}3+\f{3^p-3}4\r)\pmod{p^2}.$$

Let $p=mf+1$ be a prime, where $m,f\in\Z^+$ and $m>1$. It is interesting to determine
$\bi{rf+sf}{rf}$ mod $p^2$ in terms of parameters arising from representations of $p$ by certain binary quadratic forms
where $r$ and $s$ are positive integers with $r+s\ls m$. The reader may consult the survey [HW] and [BEW, Chapter 9]
for results and methods on this problem with $m=3,\ldots,16,20,24$.
See also [Yo] for related results obtained via the Gross-Koblitz formula.

Let $d$ be a squarefree positive integer. It is known that the
imaginary quadratic field $\Q(\sqrt{-d})$ has class number one if and only if $d\in\{1,2,3,7,11,19,43,67,163\}$.

Let $d\in\{7,11,19,43,67,163\}$. Then $K=\Q(\sqrt{-d})$ has class number one.
Let $p$ be an odd prime with $(\f{-d}p)=1$. By algebraic number theory,
$p$ splits in the ring $O_K$ of algebraic integers in $K$. As $O_K$ is a principal ideal domain,
there are unique positive integers $x$ and $y$ with $x\eq y\ (\mo\ 2)$ such that
$$p=\f{x+y\sqrt{-d}}2\times\f{x-y\sqrt{-d}}2$$
and hence $4p=x^2+dy^2$.

Let $d\in\Z^+$ be squarefree.  Suppose that $p$ is an odd prime with
$4p=x^2+dy^2$ for some $x,y\in\Z$. If both $x$ and $y$ are odd, then
$d+1\eq4p\eq4\pmod{8}$ and hence $d\eq 3\pmod8$. Thus, if
$d\not\eq3\pmod8$, then $x$ and $y$ are both even, and
$p=(x/2)^2+d(y/2)^2$.

Let $p$ be an odd prime.  In 1977, A. R. Rajwade [Ra] proved that
$$\align \sum_{x=0}^{p-1}\l(\f{x^3+21x^2+112x}p\r)
=\cases -2x(\f x7)&\t{if}\ (\f p7)=1\ \&\ p=x^2+7y^2\ (x,y\in\Z),
\\0&\t{if}\ (\f p7)=-1,\ \t{i.e.},\ p\eq3,5,6\pmod7.\endcases
\endalign$$
Using character sums, Jacobi (cf. [HW]) proved that if $p=11f+1$  with $f\in\Z^+$ and  $4p=x^2+11y^2$
with $x,y\in\Z$ and $x\eq2\pmod{11}$ then
$$x\eq\bi{6f}{3f}\bi{3f}f\bigg/\bi{4f}{2f}\pmod p.$$
In 1982 J. C. Parnami and Rajwade [PR] showed that
$$\align&\sum_{x=0}^{p-1}\l(\f{x^3-96\cdot11x+112\cdot11^2}p\r)
\\=&\cases (\f2p)(\f x{11})x&\t{if}\ (\f p{11})=1\ \&\ 4p=x^2+11y^2\ (x,y\in\Z),
\\0&\t{if}\ (\f p{11})=-1,\ \t{i.e.},\ p\eq2,6,7,8,10\pmod {11}.\endcases
\endalign$$
Via elliptic curves with complex multiplication, it is also known that (cf. [RPR], [JM] and [PV])
for $d=19,43,67,163$ we have
$$\sum_{x=0}^{p-1}\l(\f{f_d(x)}p\r)=\cases(\f 2p)(\f xd)x&\t{if}\ (\f pd)=1\ \&\ 4p=x^2+dy^2\ (x,y\in\Z),
\\0&\t{if}\ (\f pd)=-1,\endcases$$ where
$$\align f_{19}(x)=&x^3-8\cdot 19x+2\cdot 19^2,
\\f_{43}(x)=&x^3-80\cdot 43x+42\cdot 43^2,
\\f_{67}(x)=&x^3-440\cdot 67x+434\cdot 67^2,
\\f_{163}(x)=&x^3-80\cdot23\cdot29\cdot163x+14\cdot11\cdot19\cdot127\cdot163^2.
\endalign$$
See also Williams [W] for other similar results.

Let $p\eq1\pmod4$ be a prime and write $p=x^2+y^2$ with $x\eq1\pmod 4$ and $y\eq0\pmod2$.
The author [Su3, Conjecture 5.5] conjectured that
$$\sum_{k=0}^{p-1}\f{\bi{2k}k^2}{8^k}\eq\sum_{k=0}^{p-1}\f{\bi{2k}k^2}{(-16)^k}\eq\l(\f 2p\r)\sum_{k=0}^{p-1}\f{\bi{2k}k^2}{32^k}
\eq\l(\f 2p\r)\l(2x-\f p{2x}\r)\ (\mo\ p^2),$$
and this was completely proved by the author's twin brother Zhi-Hong Sun [S1] with the help of Legendre polynomials.

Let $n\in\N=\{0,1,2,\ldots\}$. In combinatorics, the central Delannoy number $D_n$ denotes
the number of lattice paths from $(0,0)$
to $(n,n)$ with steps $(1,0),(0,1)$ and $(1,1)$ (cf. [S]). It is known that
$$D_n=\sum_{k=0}^n\bi nk\bi{n+k}k=\sum_{k=0}^n\bi{n+k}{2k}\bi{2k}k=[x^n](x^2+3x+2)^n.$$
(Throughout this paper by $[x^n]P(x)$ we mean the coefficient
of $x^n$ in the power series of $P(x)$.)

Let $p=2n+1$ be an odd prime. For $k=n+1,\ldots,p-1$, we clearly have
$$\bi{2k}k=\f{(2k)!}{(k!)^2}\eq0\pmod{p}.$$
For $k=0,\ldots,n$, Z. H. Sun [S1] noted that
$$\align\bi{n+k}{2k}=&\f{\prod_{0<j\ls k}(p^2-(2j-1)^2)}{4^k(2k)!}
\\\eq&\f{\prod_{0<j\ls k}(-(2j-1)^2)}{4^k(2k)!}
=\f{\bi{2k}k}{(-16)^k}\ (\mo\ p^2).
\endalign$$
and this implies van Hamme's observation ([vH])
$$\bi nk\bi{n+k}k=\bi{n+k}{2k}\bi{2k}k\eq\f{\bi{2k}k^2}{(-16)^k}\pmod{p^2}.$$
Therefore
$$D_{(p-1)/2}=D_n\eq\sum_{k=0}^n\f{\bi{2k}k^2}{(-16)^k}\eq\sum_{k=0}^{p-1}\f{\bi{2k}k^2}{(-16)^k}\pmod{p^2}.$$
Since
$$\bi{(p-1)/2}k\eq\bi{-1/2}k=\f{\bi{2k}k}{(-4)^k}\pmod{p}\quad\t{for}\
k=0,1,\ldots,p-1,$$ we have
$$\sum_{k=0}^{p-1}\f{\bi{2k}k^2}{(-16)^k}\eq\sum_{k=0}^{n}(-1)^k\bi{n}k^2=[x^n](1+x)^n(1-x)^n=[x^n](1-x^2)^n\pmod
p.$$
Thus, if $p\eq3\pmod4$ (i.e.,
$2\nmid n$) then we have
$\sum_{k=0}^{p-1}\bi{2k}k^2/(-16)^k\eq0\pmod p$. If $p\eq1\pmod 4$
and $p=x^2+y^2$ with $x\eq1\pmod 4$ and $y\eq0\pmod 2$, then
$$\sum_{k=0}^{p-1}\f{\bi{2k}k^2}{(-16)^k}\eq(-1)^{(p-1)/4}\bi{(p-1)/2}{(p-1)/4}\eq(-1)^{(p-1)/4}2x\pmod p$$
by applying Gauss' congruence. This determines $\sum_{k=0}^{p-1}\bi{2k}k^2/(-16)^k$ mod $p$
in a simple way. Similarly, in 2009 the author determined $\sum_{k=0}^{p-1}\bi{2k}k^2/8^k$ and $\sum_{k=0}^{p-1}\bi{2k}k^2/32^k$
modulo $p$ by noting that
$$\align\sum_{k=0}^{p-1}\f{\bi{2k}k^2}{8^k}\eq&\sum_{k=0}^n 2^k\bi nk\bi n{n-k}=[x^n](1+2x)^n(1+x)^n
\\=&[y^n](y^2+3y+2)^n=D_n\eq\sum_{k=0}^{p-1}\f{\bi{2k}k^2}{(-16)^k}\pmod p
\endalign$$ and
$$\align\sum_{k=0}^{p-1}\f{\bi{2k}k^2}{32^k}\eq&\sum_{k=0}^n\f1{2^k}\bi nk^2=[x^n]\l(1+\f x2\r)^n(1+x)^n
\\=&2^{-n}[x^n](x^2+3x+2)^n=\f{D_n}{2^n}\eq\l(\f 2p\r)\sum_{k=0}^{p-1}\f{\bi{2k}k^2}{(-16)^k}\pmod{p}.
\endalign$$

Recently the author established the following theorem.
\proclaim{Theorem 1.1} {\rm (Z. W. Sun [Su5])} Let $p\eq1\pmod4$ be a prime and write $p=x^2+y^2$ with $x\eq1\pmod4$ and $y\eq0\pmod2$.
Then we can determine $x$ mod $p^2$ in the following way:
$$(-1)^{(p-1)/4}\,x\eq\sum_{k=0}^{(p-1)/2}\f{k+1}{8^k}\bi{2k}k^2
\eq\sum_{k=0}^{(p-1)/2}\f{2k+1}{(-16)^k}\bi{2k}k^2\pmod{p^2}.\tag1.1$$
\endproclaim

Concerning the representations $p=x^2+3y^2$ and $p=x^2+7y^2$, we have formulated a conjecture similar to Theorem 1.1.
\proclaim{Conjecture 1.1} {\rm (Sun [Su4])} Let $p$ be an odd prime.

{\rm (i)} If $p\eq1\pmod3$ and $p=x^2+3y^2$ with $x,y\in\Z$ and $x\eq1\pmod3$, then we can determine $x$ mod $p^2$
in the following way:
$$\sum_{k=0}^{p-1}\f{k+1}{48^k}\bi{2k}k\bi{4k}{2k}\eq x\ (\mo\ p^2).\tag1.2$$

{\rm (ii)} If $p\eq7\ (\mo\ 12)$ and $p=x^2+3y^2$ with $y\eq1\ (\mo\ 4)$, then
we can determine $y$ mod $p^2$ via the congruence
$$\sum_{k=0}^{p-1}\l(\f k3\r)\f{k\bi{2k}k^2}{(-16)^k}\eq (-1)^{(p+1)/4}y\ (\mo\ p^2).\tag1.3$$

{\rm (iii)} If $(\f p7)=1$ and $p=x^2+7y^2$ with $x,y\in\Z$ and $(\f x7)=1$, then we can determine $x$ mod $p^2$
in the following way:
$$\sum_{k=0}^{p-1}\f{k+8}{63^k}\bi{2k}k\bi{4k}{2k}\eq 8\l(\f p3\r)x\ (\mo\ p^2).\tag1.4$$
\endproclaim

For congruences concerning the sum $\sum_{k=0}^{p-1}\bi{2k}k/m^k$ modulo powers of a prime $p$
(where $m\in\Z$ and $m\not\eq0\pmod p$), the reader may consult [PS], [ST1], [ST2], [Su1] and [Su2].

Let $d\in\Z^+$ be squarefree and let $p$ be an odd prime with $(\f{-d}p)=1$.
If $p$ or $4p$ can be written in the form $x^2+dy^2$ (in such a case we always assume $x,y\in\Z$ even though
we may not mention it to avoid lines overfull),
how to determine $x^2$ modulo $p^2$ in a simple and explicit form?
We will focus on this question in this survey and present various conjectures and results on the author's following problems.

\proclaim{Problem 1.1} Given a squarefree positive integer $d$,
find integers $a_0,a_1,\ldots$ such that for sufficiently large primes $p$ we have
$$\sum_{k=0}^{p-1}a_k\eq\cases x^2-2p\pmod {p^2}&\t{if}\ 4p=x^2+dy^2\ (\t{and}\ 4\nmid x\ \t{if}\ d=1),
\\0\pmod{p^2}&\t{if}\ (\f{-d}p)=-1.\endcases$$
\endproclaim

If one thinks that the integral condition of $a_0,a_1,a_2,\ldots$ in Problem 1.1 is too strict,
we may study the following easier problem.

\proclaim{Problem 1.2} Given a squarefree $d\in\Z^+$,
find rational numbers $a_0,a_1,a_2,\ldots$ with denominators not divisible by large primes
such that for sufficiently large primes $p$ we have
$$\sum_{k=0}^{p-1}a_k\eq\cases x^2-2p\pmod {p^2}&\t{if}\ 4p=x^2+dy^2\ (\t{and}\ 4\nmid x\ \t{if}\ d=1),
\\0\pmod{p^2}&\t{if}\ (\f{-d}p)=-1.\endcases$$
\endproclaim

We find that Problems 1.1 and 1.2 have affirmative answers
for most of those $d\in\Z^+$ with the imaginary quadratic field $\Q(\sqrt {-d})$ having class number 1 or 2 or 4.

 Now we give suggested solutions to the problems with $d=7,15,427$ as examples.

\proclaim{Conjecture 1.2 {\rm (Sun [Su3])}} Let $p$ be an odd prime. Then
$$\sum_{k=0}^{p-1}\bi{2k}k^3
\eq\cases4x^2-2p\ (\mo\ p^2)&\t{if}\ (\f p7)=1\ \&\ p=x^2+7y^2\ (x,y\in\Z),
\\0\ (\mo\ p^2)&\t{if}\ (\f p7)=-1,\ \t{i.e.,}\ p\eq3,5,6\ (\mo\ 7).\endcases
\tag1.5$$
\endproclaim
\Remark\ 1.1. The congruence modulo $p$ can be easily deduced from (6) in Ahlgren [A, Theorem 5].
Recently the author's twin brother Z. H. Sun [S2] confirmed the conjecture in the case $(\f p7)=-1$
via Legendre polynomials.

\medskip
The imaginary quadratic fields $\Q(\sqrt{-15})$ and $\Q(\sqrt{-427})$ have class number two.
Our next two conjectures
provide affirmative answers to Problem 1.1 with $d=15$ and Problem 1.2 with $d=427$.

\proclaim{Conjecture 1.3} For $k=0,1,2,\ldots$ let $T_k$ denote the trinomial coefficient
defined as the coefficient of $x^k$ in $(x^2+x+1)^k$.
For any prime $p>3$, we have
$$\aligned&\sum_{k=0}^{p-1}(-1)^k\bi{2k}k^2T_k
\\\eq&\cases 4x^2-2p\pmod{p^2}&\t{if}\ p\eq1,4\pmod{15}\ \&\ p=x^2+15y^2\ (x,y\in\Z),
\\2p-12x^2\pmod{p^2}&\t{if}\ p\eq2,8\pmod{15}\ \&\ p=3x^2+5y^2\ (x,y\in\Z),
\\0\pmod{p^2}&\t{if}\ (\f p{15})=-1,\ \t{i.e.},\ p\eq 7,11,13,14\pmod{15},\endcases
\endaligned$$
and
$$\sum_{k=0}^{p-1}(105k+44)(-1)^k\bi{2k}k^2T_k\eq p\l(20+24\l(\f p3\r)(2-3^{p-1})\r)\pmod{p^3}.$$
Also,
$$\f1{2n\bi{2n}n}\sum_{k=0}^{n-1}(-1)^{n-1-k}(105k+44)\bi{2k}k^2T_k\in\Z^+\quad\t{for all}\ n=1,2,3,\ldots.$$
\endproclaim

\proclaim{Conjecture 1.4} Let $p$ be an odd prime with $p\not=11$.
Then
$$\align&\sum_{n=0}^{p-1}\f1{(-22)^{3n}}\sum_{k=0}^n\bi nk^2\bi{n+k}k\bi{2k}n
\\\eq&\cases x^2-2p\pmod{p^2}&\t{if}\ (\f p7)=(\f p{61})=1\ \&\ 4p=x^2+427y^2\ (x,y\in\Z),
\\2p-7x^2\pmod{p^2}&\t{if}\ (\f p7)=(\f p{61})=-1\ \&\ 4p=7x^2+61y^2\ (x,y\in\Z),
\\0\pmod{p^2}&\t{if}\ (\f p{427})=-1.\endcases
\endalign$$
We also have
$$\sum_{n=0}^{p-1}\f{11895n+1286}{(-22)^{3n}}
\sum_{k=0}^n\bi nk^2\bi{n+k}k\bi{2k}n\eq 1286p\l(\f p7\r)\pmod{p^2}.$$
\endproclaim

L. van Hamme [vH] observed that certain hypergeometric series have their $p$-adic analogues, e.g., motivated by the
identity (of Ramanujan's type) $$\sum_{k=0}^\infty\f{4k+1}{(-64)^k}\bi{2k}k^3=\f 2{\pi}$$
he conjectured that
$$\sum_{k=0}^{p-1}\f{4k+1}{(-64)^k}\bi{2k}k^3\eq p\l(\f{-1}p\r)\pmod{p^3}$$
for any odd prime $p$ (and this was confirmed by E. Mortenson [M]).
The reader may consult [BBC], [Be], [BoBo] and [ChCh] for more Ramanujan-type series for $1/\pi$.

The author has the following viewpoint (or hypothesis) the initial version of which appeared in
his message to Number Theory Mailing List (sent on March 30, 2010) available from the website
{\tt http://listserv.nodak.edu/cgi-bin/wa.exe?}
{\tt A2=ind1003$\&$L=nmbrthry$\&$T=0$\&$P=1668}.

\proclaim{Philosophy about Series for $1/\pi$} {\rm (i)} Given a `regular' identity of the form
$$\sum_{k=0}^\infty(bk+c)\f{a_k}{m^k}=\f C{\pi},$$
where $a_k,b,c,m\in\Z,$ $bm$ is nonzero and $C^2$ is rational, we must have
$$\sum_{k=0}^{n-1}(bk+c)a_km^{n-1-k}\eq0\pmod{n}$$
for any positive integer $n$. Furthermore,
there exist an integer $m'$ (often $m'=m$)
and a squarefree positive integer $d$
with the class number of $\Q(\sqrt{-d})$ in $\{1,2,2^2,\ldots\}$ (and with $C/\sqrt d$ often rational)
such that either $d>1$ and for any prime $p>3$ not dividing $dm$ we have
$$\sum_{k=0}^{p-1}\f{a_k}{m^k}
\eq\cases (\f {m'}p)(x^2-2p)\pmod {p^2}&\t{if}\ 4p=x^2+dy^2\ (x,y\in\Z),
\\0\pmod{p^2}&\t{if}\ (\f{-d}p)=-1,\endcases$$
or $d=1$, $\gcd(15,m)>1$, and $\sum_{k=0}^{p-1}a_k/{m^k}\eq0\pmod{p^2}$
for any prime $p\eq3\pmod4$ with $p\nmid 3m$.

{\rm (ii)} Let $b,c,m,a_0,a_1,\ldots$ be integers with $bm$ nonzero and the series $\sum_{k=0}^\infty(bk+c)a_k/m^k$ convergent.
Suppose that there are $d\in\Z^+$, $d'\in\Z$, and rational numbers $c_0$ and $c_1$ such that
$$\sum_{k=0}^{p-1}(bk+c)\f{a_k}{m^k}\eq p\l(c_0\l(\f{-d}p\r)+c_1\l(\f{d'}p\r)\r)\pmod{p^2}$$
for all sufficiently large primes $p$. If $d'\gs0$, then
$$\sum_{k=0}^\infty(bk+c)\f{a_k}{m^k}=\f C{\pi}$$
for some $C$ with $C^2$ rational (and with $C/\sqrt d$ rational if $c_0\not=0$).
If $d'=-d_1<0$, then there are rational numbers $\lambda_0$ and $\lambda_1$ such that
$$\sum_{k=0}^\infty(bk+c)\f{a_k}{m^k}=\f{\lambda_0\sqrt d+\lambda_1\sqrt{d_1}}{\pi}.$$
\endproclaim
\Remark\ 1.2. We don't have a precise description of the meaning of the word {\it regular}
in the above philosophy, but almost all identities of the stated form are regular.
\medskip

Now we illustrate the philosophy by the author's next three conjectures.

\proclaim{Conjecture 1.5} For any prime $p>3$, we have
$$\align&\sum_{n=0}^{p-1}\f{\bi{2n}n}{256^n}\sum_{k=0}^n\bi{2k}k^2\bi{2(n-k)}{n-k}12^{n-k}
\\\eq&\cases4x^2-2p\pmod{p^2}&\t{if}\ p\eq1,7\pmod{24}\ \&\ p=x^2+6y^2,
\\8x^2-2p\pmod{p^2}&\t{if}\ p\eq5,11\pmod{24}\ \&\ p=2x^2+3y^2,
\\0\pmod{p^2}&\t{if}\ (\f{-6}p)=-1,\endcases
\endalign$$
and
$$\sum_{n=0}^{p-1}\f{6n-1}{256^n}\bi{2n}n\sum_{k=0}^n\bi{2k}k^2\bi{2(n-k)}{n-k}12^{n-k}\eq -p\pmod{p^2}.$$
Also,
$$\sum_{n=0}^\infty\f{6n-1}{256^n}\bi{2n}n\sum_{k=0}^n\bi{2k}k^2\bi{2(n-k)}{n-k}12^{n-k}=\f{8{\sqrt3}}{\pi}.$$
\endproclaim

For $b,c\in\Z$ the $n$th generalized central trinomial coefficient $T_n(b,c)$ is defined as $[x^n](x^2+bx+c)^n$.
In [Su8] the author conjectured that
$$\sum_{k=0}^\infty\f{80k+13}{(-168^2)^{k}}\bi{4k}{2k}\bi{2k}kT_{k}(7,4096)=\f{14\sqrt{210}+21\sqrt{42}}{8\pi}.$$
Here are related conjectural congruences.

\proclaim{Conjecture 1.6}  Let $p>3$ be a prime with $p\not=7$. Then
$$\align&\l(\f{-42}p\r)\sum_{k=0}^{p-1}\f{\bi{4k}{2k}\bi{2k}kT_{k}(7,4096)}{(-168^2)^{k}}
\\\eq&\cases 4x^2-2p\pmod{p^2}&\t{if}\ p\eq1,4\pmod{15}\ \&\
p=x^2+15y^2,
\\12x^2-2p\pmod{p^2}&\t{if}\ p\eq2,8\pmod{15}\ \&\ p=3x^2+5y^2,
\\0\pmod{p^2}&\t{if}\ (\f p{15})=-1.\endcases
\endalign$$
And
$$\sum_{k=0}^{p-1}\f{80k+13}{(-168^2)^k}\bi{4k}{2k}\bi{2k}kT_{k}(7,4096)
\eq p\l(3\l(\f{{-42}}p\r)+10\l(\f{{-210}}p\r)\r)\pmod{p^2}.$$
\endproclaim

The imaginary quadratic field $\Q(\sqrt{-105})$ has class number eight, so our following conjecture is of particular interest.

\proclaim{Conjecture 1.7} For any prime $p>5$, we have
$$\align&\sum_{n=0}^{p-1}\f{\bi{2n}n}{2160^n}\sum_{k=0}^n\bi nk\bi{n+2k}{2k}\bi{2k}k(-324)^{n-k}
\\\eq&\cases 4x^2-2p\,(\mo\ p^2)&\t{if}\, (\f{-1}p)=(\f p3)=(\f p5)=(\f p7)=1,\, p=x^2+105y^2,
\\2x^2-2p\,(\mo\ p^2)&\t{if}\, (\f{-1}p)=(\f p7)=1,\,(\f p3)=(\f p5)=-1,\, 2p=x^2+105y^2,
\\2p-12x^2\,(\mo\ p^2)&\t{if}\,(\f{-1}p)=(\f p3)=(\f p5)=(\f p7)=-1,\, p=3x^2+35y^2,
\\2p-6x^2\,(\mo\ p^2)&\t{if}\, (\f{-1}p)=(\f p7)=-1,\, (\f p3)=(\f p5)=1,\, 2p=3x^2+35y^2,
\\20x^2-2p\,(\mo\ p^2)&\t{if}\, (\f{-1}p)=(\f p5)=1,\,(\f p3)=(\f p7)=-1,\, p=5x^2+21y^2,
\\10x^2-2p\,(\mo\ p^2)&\t{if}\, (\f{-1}p)=(\f p3)=1,\, (\f p5)=(\f p7)=-1,\, 2p=5x^2+21y^2,
\\28x^2-2p\,(\mo\ p^2)&\t{if}\, (\f{-1}p)=(\f p5)=-1,\,(\f p3)=(\f p7)=1,\, p=7x^2+15y^2,
\\14x^2-2p\,(\mo\ p^2)&\t{if}\, (\f{-1}p)=(\f p3)=-1,\, (\f p5)=(\f p7)=1,\, 2p=7x^2+15y^2,
\\0\pmod{p^2}&\t{if}\ (\f{-105}p)=-1,
\endcases
\endalign$$
and
$$\align&\sum_{n=0}^{p-1}\f{357n+103}{2160^n}\bi{2n}n\sum_{k=0}^n\bi nk\bi{n+2k}{2k}\bi{2k}k(-324)^{n-k}
\\\quad\qquad&\eq p\l(54\l(\f{-1}p\r)+49\l(\f {15}p\r)\r)\pmod{p^2}.
\endalign$$
Also, $$\sum_{n=0}^\infty\f{357n+103}{2160^n}\bi{2n}n\sum_{k=0}^n\bi nk\bi{n+2k}{2k}\bi{2k}k(-324)^{n-k}=\f{90}{\pi}.$$
\endproclaim

We mention that series for $1/\pi$ are usually difficult to prove; most proofs of them involve
modular functions or elliptic integrals.

From Section 2 we will give various conjectures and results on Problems 1.1 and 1.2.
Problem 1.1 for $d=1$ already has a positive answer. We suggest positive answers to Problem 1.1 for
$$d\in\{2,3,5,6,7,10,13,15,22,30,37,58,70,85,130,190\}.$$
We also formulate many conjectures concerning Problem 1.2; in
particular, we give explicit conjectural positive answers for those
squarefree positive integers $d$ with $\Q(\sqrt{-d})$ having class
number at most two except for $d=187,403$. Note that $\Q(\sqrt{-d})$
has class number two if and only if
$$d\in\{5, 6, 10, 13, 15, 22, 35, 37, 51, 58, 91, 115, 123, 187, 235, 267, 403, 427\}.$$
For each of
$$d\in\{21,30,33,42,57,70,78,85,93,102,130,133,177,190\},$$
the quadratic field $\Q(\sqrt{-d})$ has class number four; for these values of $d$
we also provide explicit conjectural positive answers to
Problem 1.2 or even Problem 1.1.

\heading{2. Using Ap\'ery polynomials and products of three binomial coefficients}
\endheading

The well-known Ap\'ery numbers (introduced by Ap\'ery in his proof
of the irrationality of $\zeta(3)=\sum_{n=1}^\infty 1/n^3$ (see [Ap] and [P]), are given by
$$A_n=\sum_{k=0}^n\bi
nk^2\bi{n+k}k^2=\sum_{k=0}^n\bi{n+k}{2k}^2\bi{2k}k^2\
(n=0,1,2,\ldots),$$
They also have close connections to modular forms (cf. K. Ono [O]).
As in [Su7]  we define Ap\'ery polynomials by
$$A_n(x)=\sum_{k=0}^n\bi nk^2\bi{n+k}k^2x^k=\sum_{k=0}^n\bi{n+k}{2k}^2\bi{2k}k^2x^k\ \ (n=0,1,2,\ldots).$$
Clearly $A_n(1)=A_n$. The author [Su7] proved that for all $n\in\Z^+$ and $x\in\Z$ we have
$$\sum_{k=0}^{n-1}(2k+1)A_k(x)\eq0\pmod{n}.$$
He also conjectured that $n\mid\sum_{k=0}^{n-1}(2k+1)(-1)^kA_k(x)$ for all $n\in\Z^+$ and $x\in\Z$, which was later confirmed by
V. J. W. Guo and J. Zeng [GZ].

The following theorem solves Problem 1.1 for $d=1$.
\proclaim{Theorem 2.1 {\rm (Sun [Su7])}} Let $p$ be an odd prime. Then
$$\aligned&\sum_{k=0}^{p-1}(-1)^kA_k(-2)\eq\sum_{k=0}^{p-1}(-1)^kA_k\l(\f 14\r)
\\\eq&\cases 4x^2-2p\ (\mo\ p^2)&\t{if}\ p\eq1\ (\mo\ 4)\ \&\ p=x^2+y^2\ (2\nmid x),
\\0\ (\mo\ p^2)&\t{if}\ p\eq3\ (\mo\ 4).\endcases
\endaligned\tag2.1$$
\endproclaim
\Remark\ 2.1. A lemma for the proof of Theorem 2.1 states that for any odd prime $p$ we have
$$\sum_{k=0}^{p-1}\f{\bi{2k}k^3}{(-8)^k}\eq\cases 4x^2-2p\ (\mo\ p^2)&\t{if}\ p\eq1\ (\mo\ 4)\ \&\ p=x^2+y^2\ (2\nmid x),
\\0\ (\mo\ p^2)&\t{if}\ p\eq3\ (\mo\ 4).\endcases$$
This was first conjectured by the author [Su4] and later confirmed by his twin brother Z.-H. Sun [S2].

\medskip

Concerning Problem 1.1 for $d=2,3$ the author made the following conjecture.

\proclaim{Conjecture 2.1 {\rm (Sun [Su7])}} Let $p$ be an odd prime. Then
$$\sum_{k=0}^{p-1}A_k\eq \cases 4x^2-2p\pmod{p^2}&\t{if}\ p\eq1,3\pmod8\ \&\ p=x^2+2y^2,
\\0\pmod{p^2}&\t{if}\ (\f{-2}p)=-1,\ \t{i.e.,}\ p\eq5,7\pmod8.\endcases\tag2.2$$
And
$$\sum_{k=0}^{p-1}(-1)^kA_k
\eq\cases4x^2-2p\ (\mo\ p^2)&\t{if}\ p\eq1\ (\mo\ 3)\ \t{and}\ p=x^2+3y^2,
\\0\ (\mo\ p^2)&\t{if}\ p\eq2\ (\mo\ 3).\endcases
\tag2.3$$
\endproclaim
\Remark\ 2.2. The author [Su7] proved the mod $p$ version of (2.2) and (2.3), and (2.3) in the case $p\eq2\pmod3$.
\medskip

The following theorem relates sums of Ap\'ery polynomials to sums of products of three binomial coefficients.
Obviously, for each $k\in\N$ we have
$$\bi{4k}{k,k,k,k}=\f{(4k)!}{(k!)^4}=\bi{2k}k^2\bi{4k}{2k}.$$

\proclaim{Theorem 2.2 {\rm (Sun [Su7])}} Let $p$ be an odd prime. Then
$$\sum_{k=0}^{p-1}(-1)^kA_k(x)\eq\sum_{k=0}^{p-1}\f{\bi{2k}k^3}{16^k}x^k\pmod{p^2}.\tag2.4$$
Also, for any $p$-adic integer $x\not\eq0\pmod p$ we have
$$\aligned\sum_{k=0}^{p-1}A_k(x)
\eq\l(\f xp\r)\sum_{k=0}^{p-1}\f{\bi{4k}{k,k,k,k}}{(256x)^k}\pmod{p}.
\endaligned\tag2.5$$
\endproclaim

Recall that Conjecture 1.2 determines $\sum_{k=0}^{p-1}\bi{2k}k^3$ mod $p^2$ for any odd prime $p$.
By Theorem 2.2, we have
$$\sum_{k=0}^{p-1}(-1)^kA_k(16)\eq \sum_{k=0}^{p-1}\bi{2k}k^3\pmod{p^2}.$$

\proclaim{Conjecture 2.2} Let $p\not=2,5$ be a prime. Then
 $$\align\sum_{k=0}^{p-1}A_k(-4)
 \eq\cases
 4x^2-2p\ (\mo\ p^2)&\t{if}\ p\eq1,9\ (\mo\ 20)\ \&\ p=x^2+5y^2,
 \\2x^2-2p\ (\mo\ p^2)&\t{if}\ p\eq3,7\ (\mo\ 20)\ \&\ 2p=x^2+5y^2,
 \\0\ (\mo\ p^2)&\t{if}\ p\eq 11,13,17,19\,(\mo\ 20).\endcases
 \endalign$$
 \endproclaim
\Remark\ 2.3. Let $p\not=2,5$ be a prime. By the theory of binary quadratic
forms (see, e.g., [C]), if $p\eq1,9\ (\mo\ 20)$ then $p=x^2+5y^2$
for some $x,y\in\Z$; if $p\eq 3,7\ (\mo\ 20)$ then $2p=x^2+5y^2$ for
some $x,y\in\Z$.

\proclaim{Conjecture 2.3} Let $p>3$ be a prime. Then
$$\sum_{k=0}^{p-1}A_k(9)
 \eq\cases
 4x^2-2p\ (\mo\ p^2)&\t{if}\ p\eq1,7\ (\mo\ 24)\ \&\ p=x^2+6y^2,
 \\2p-8x^2\ (\mo\ p^2)&\t{if}\ p\eq5,11\ (\mo\ 24)\ \&\ p=2x^2+3y^2,
 \\0\ (\mo\ p^2)&\t{if}\ p\eq 13,17,19,23\, (\mo\ 24).
 \endcases$$
 \endproclaim

\proclaim{Conjecture 2.4} Let $p>3$ be a prime with $p\not=11$. Then
$$\align&\sum_{k=0}^{p-1}A_k(99^2)\eq\sum_{k=0}^{p-1}\f{\bi{4k}{k,k,k,k}}{1584^{2k}}
\\\eq&\cases 4x^2-2p\ (\mo\ p^2)&\t{if}\ (\f{-11}p)=(\f{2}p)=1\ \&\ p=x^2+22y^2,
\\2p-8x^2\ (\mo\ p^2)&\t{if}\ (\f{-11}p)=(\f{2}p)=-1\ \&\ p=2x^2+11y^2,
\\0\ (\mo\ p^2)&\t{if}\ (\f{-11}p)=-(\f{2}p).\endcases
\endalign$$
And
$$\align&\sum_{k=0}^{p-1}A_k(99^4)\eq\sum_{k=0}^{p-1}\f{\bi{4k}{k,k,k,k}}{396^{4k}}
\\\eq&\cases 4x^2-2p\ (\mo\ p^2)&\t{if}\ (\f{29}p)=(\f{-2}p)=1\ \&\ p=x^2+58y^2,
\\2p-8x^2\ (\mo\ p^2)&\t{if}\ (\f{29}p)=(\f{-2}p)=-1\ \&\ p=2x^2+29y^2,
\\0\ (\mo\ p^2)&\t{if}\ (\f{-58}p)=-1.\endcases
\endalign$$
Furthermore, for $n=2,3,4,\ldots$ we have
$$a_n:=\f1{2n(2n+1)\bi{2n}n}\sum_{k=0}^{n-1}(280k+19)\bi{4k}{k,k,k,k}1584^{2(n-1-k)}\in\Z$$
unless $2n+1$ is a power of $3$ in which case $3a_n\in\Z\sm3\Z$.
Also, for $n=2,3,4,\ldots$ we have
$$b_n:=\f1{2n(2n+1)\bi{2n}n}\sum_{k=0}^{n-1}(26390k+1103)\bi{4k}{k,k,k,k}396^{4(n-1-k)}\in\Z$$
unless $2n+1$ is a power of $3$ in which case $3b_n\in\Z\sm3\Z$.
\endproclaim
\Remark\ 2.4. Ramanujan (cf. [Be, p.\,354]) found that
$$\sum_{k=0}^\infty\f{280k+19}{1584^{2k}}\bi{4k}{k,k,k,k}=\f{2\times99^2}{\pi\sqrt{11}}$$
and
$$\sum_{k=0}^\infty\f{26390k+1103}{396^{4k}}\bi{4k}{k,k,k,k}=\f{99^2}{2\pi\sqrt2}.$$
The reader may consult [BB], [BBC] and [CC] for other known Ramanujan-type series not mentioned in this paper.

\proclaim{Conjecture 2.5} Let $p>3$ be a prime.
Then
$$\align&\sum_{k=0}^{p-1}A_k(-324)\eq
\l(\f{-1}p\r)\sum_{k=0}^{p-1}\f{\bi{4k}{k,k,k,k}}{(-2^{10} 3^4)^k}
\\\eq&\cases 4x^2-2p\ (\mo\ p^2)&\t{if}\ (\f{13}p)=(\f{-1}p)=1\ \&\ p=x^2+13y^2,
\\2x^2-2p\ (\mo\ p^2)&\t{if}\ (\f{13}p)=(\f{-1}p)=-1\ \&\ 2p=x^2+13y^2,
\\0\ (\mo\ p^2)&\t{if}\ (\f{13}p)=-(\f{-1}p).\endcases
\endalign$$
We also have
$$\sum_{k=0}^{p-1}\f{260k+23}{(-2^{10}3^4)^k}\bi{4k}{k,k,k,k}\eq 23p\l(\f{-1}p\r)+\f 53 p^3E_{p-3}\ (\mo\ p^4),$$
where $E_0,E_1,E_2,\ldots$ are Euler numbers.
Furthermore, for $n=2,3,4,\ldots$ we have
$$c_n:=\f1{2n(2n+1)\bi{2n}n}\sum_{k=0}^{n-1}(260k+23)\bi{4k}{k,k,k,k}(-82944)^{n-1-k}\in\Z$$
unless $2n+1$ is a power of $3$ in which case $3c_n\in\Z\sm3\Z$.
\endproclaim

\proclaim{Conjecture 2.6} Let $p>3$ be a prime with $p\not=7$. Then
$$\align&\sum_{k=0}^{p-1}A_k(-882^2)\eq\l(\f{-1}p\r)\sum_{k=0}^{p-1}\f{\bi{4k}{k,k,k,k}}{(-2^{10} 21^4)^k}
\\\eq&\cases 4x^2-2p\ (\mo\ p^2)&\t{if}\ (\f{37}p)=(\f{-1}p)=1\ \&\ p=x^2+37y^2,
\\2x^2-2p\ (\mo\ p^2)&\t{if}\ (\f{37}p)=(\f{-1}p)=-1\ \&\ 2p=x^2+37y^2,
\\0\ (\mo\ p^2)&\t{if}\ (\f{-37}p)=-1.\endcases
\endalign$$
Furthermore, for $n=2,3,4,\ldots$ we have
$$d_n:=\f1{2n(2n+1)\bi{2n}n}\sum_{k=0}^{n-1}(21460k+1123)\bi{4k}{k,k,k,k}(-2^{10}21^4)^{n-1-k}\in\Z$$
unless $2n+1$ is a power of $3$ in which case $3d_n\in\Z\sm3\Z$.
\endproclaim

We will not list our similar conjectures on $\sum_{k=0}^{p-1}A_k(x)$ mod $p^2$ with $x$ among
$$-48,\ 81,\ 2401,\ -25920,\ -\f 9{16},\ \f{81}{32},\ \f{81}{256},\ -\f{3969}{256}.$$

Recall that $\Q(\sqrt{-11})$ has class number one. The author's following conjecture provides a positive answer to Problem
1.2 with $d=11$.

\proclaim{Conjecture 2.7 {\rm (Sun [Su3])}} Let $p$ be an odd prime.
 Then
 $$\sum_{k=0}^{p-1}\f{\bi{2k}k^2\bi{3k}k}{64^k}
 \eq\cases
 x^2-2p\ (\mo\ p^2)&\t{if}\ (\f p{11})=1\ \&\ 4p=x^2+11y^2,
 \\0\ (\mo\ p^2)&\t{if}\ (\f p{11})=-1.
 \endcases$$
\endproclaim
\Remark\ 2.5.  Recently Z. H. Sun [S3] confirmed the mod $p$ version of the congruence.

Note that [Su4] contains some of the author's conjectural congruences related to Ramanujan-type series.
Below we present more such conjectures.

\proclaim{Conjecture 2.8}  Let $p>3$ be a prime. Then
$$\sum_{k=0}^{p-1}\f{\bi{6k}{3k}\bi{3k}{k,k,k}}{(-96)^{3k}}
\eq\cases(\f{-6}p)(x^2-2p)\ (\mo\ p^2)&\t{if}\ (\f p{19})=1\ \&\ 4p=x^2+19y^2,
\\0\ (\mo\ p^2)&\t{if}\ (\f p{19})=-1.\endcases$$
Also,
$$\sum_{k=0}^{p-1}\f{342k+25}{(-96)^{3k}}\bi{6k}{3k}\bi{3k}{k,k,k}\eq 25p\l(\f{-6}p\r)\ (\mo\ p^3).$$
\endproclaim
\Remark\ 2.6. The conjectural congruences modulo $p$ have been confirmed by Z. H. Sun [S3].
D. V. Chudnovsky and G. V. Chunovsky [CC]
obtained that
$$\sum_{k=0}^\infty\f{342k+25}{(-96)^{3k}}\bi{6k}{3k}\bi{3k}{k,k,k}=\f{32\sqrt6}{\pi}.$$

\proclaim{Conjecture 2.9} Let $p>5$ be a prime. Then
$$\sum_{k=0}^{p-1}\f{\bi{6k}{3k}\bi{3k}{k,k,k}}{(-960)^{3k}}
\eq\cases(\f{p}{15})(x^2-2p)\ (\mo\ p^2)&\t{if}\ (\f p{43})=1\ \&\ 4p=x^2+43y^2,
\\0\ (\mo\ p^2)&\t{if}\ (\f p{43})=-1.\endcases$$
If $p\not=11$, then
$$\sum_{k=0}^{p-1}\f{\bi{6k}{3k}\bi{3k}{k,k,k}}{(-5280)^{3k}}
\eq\cases(\f{-330}{p})(x^2-2p)\ (\mo\ p^2)&\t{if}\ (\f p{67})=1\ \&\ 4p=x^2+67y^2,
\\0\ (\mo\ p^2)&\t{if}\ (\f p{67})=-1.\endcases$$
\endproclaim
\Remark\ 2.7. The conjectural congruences modulo $p$ has been confirmed by Z. H. Sun [S3].

\proclaim{Conjecture 2.10} Let $p>5$ be a prime with $p\not=23,29$. Then
$$\align&\sum_{k=0}^{p-1}\f{\bi{6k}{3k}\bi{3k}{k,k,k}}{(-640320)^{3k}}
\\\eq&\cases(\f{-10005}{p})(x^2-2p)\ (\mo\ p^2)&\t{if}\ (\f p{163})=1\ \&\ 4p=x^2+163y^2,
\\0\ (\mo\ p^2)&\t{if}\ (\f p{163})=-1.\endcases
\endalign$$
Also,
$$\sum_{k=0}^{p-1}\f{545140134k+13591409}{(-640320)^{3k}}\bi{6k}{3k}\bi{3k}{k,k,k}\eq 13591409p\l(\f{-10005}p\r)\ (\mo\ p^3).$$
\endproclaim
\Remark\ 2.8. The conjectural congruence modulo $p$ has been confirmed by Z. H. Sun [S3].
D. V. Chudnovsky and G. V. Chudnovsky [CC]
got the formula
$$\sum_{k=0}^\infty\f{545140134k+13591409}{(-640320)^{3k}}\bi{6k}{3k}\bi{3k}{k,k,k}=\f{3\times53360^2}{2\pi\sqrt{10005}},$$
which enabled them to hold the record for the calculation of $\pi$ during 1989-1994.
\medskip

\proclaim{Conjecture 2.11} Let $p>3$ be a prime with $p\not=17$. Then
$$\align&\sum_{n=0}^{p-1}\f{\bi{2n}n}{68^n}\sum_{k=0}^n\bi nk^364^k\eq\sum_{k=0}^{p-1}\f{\bi{2k}{k}^2\bi{3k}{k}}{(-12^3)^k}
\\\eq&\cases x^2-2p\pmod{p^2}&\t{if}\ (\f p3)=(\f p{17})=1\ \&\ 4p=x^2+51y^2,
\\2p-3x^2\pmod{p^2}&\t{if}\ (\f p3)=(\f p{17})=-1\ \&\ 4p=3x^2+17y^2,
\\0\pmod{p^2}&\t{if}\ (\f p{51})=-1.\endcases
\endalign$$
Also,
$$\sum_{n=0}^{p-1}\f{36n+19}{68^n}\bi{2n}n\sum_{k=0}^n\bi nk^364^k\eq 19p\l(\f p{17}\r)\pmod{p^2},$$
and $$\f1{p^a}\sum_{k=0}^{p^a-1}\f{51k+7}{(-12^3)^k}\bi{2k}{k}^2\bi{3k}{k}
\eq 7\l(\f{p^a}3\r)+\l(\f{p^{a-1}}3\r)\f56p^2B_{p-2}\l(\f13\r)\ (\mo\ p^3)$$
for every $a\in\Z^+$, where $B_n(x)$ denotes the Bernoulli polynomial of degree $n$.
\endproclaim
\Remark\ 2.9. Ramanujan [R] found that
$$\sum_{k=0}^\infty\f{51k+7}{(-12^3)^{k}}\bi{2k}{k}^2\bi{3k}{k}=\f{12\sqrt3}{\pi}.$$

\proclaim{Conjecture 2.12} Let $p>3$ be a prime.

{\rm (i)} We have
$$\align\sum_{k=0}^{p-1}\f{\bi{2k}k^2\bi{3k}k}{(-48)^{3k}}
\eq\cases x^2-2p\ (\mo\ p^2)&\t{if}\ (\f{p}3)=(\f{p}{41})=1\ \&\ 4p=x^2+123y^2,
\\2p-3x^2\ (\mo\ p^2)&\t{if}\ (\f{p}3)=(\f{p}{41})=-1\ \&\ 4p=3x^2+41y^2,
\\0\ (\mo\ p^2)&\t{if}\ (\f{p}{123})=-1.\endcases
\endalign$$
Also,
$$\f1{p^a}\sum_{k=0}^{p^a-1}\f{615k+53}{(-48)^{3k}}\bi{2k}k^2\bi{3k}k\eq 53\l(\f {p^a}3\r)
+\l(\f{p^{a-1}}3\r)\f5{12}p^2B_{p-2}\l(\f13\r)\ (\mo\ p^3)$$
for any positive integer $a$.

{\rm (ii)} Suppose $p>5$. Then
$$\align&\sum_{k=0}^{p-1}\f{\bi{2k}k^2\bi{3k}k}{(-300)^{3k}}
\eq\cases x^2-2p\ (\mo\ p^2)&\t{if}\ (\f{p}3)=(\f{p}{89})=1\ \&\ 4p=x^2+267y^2,
\\2p-3x^2\ (\mo\ p^2)&\t{if}\ (\f{p}3)=(\f{p}{89})=-1\ \&\ 4p=3x^2+89y^2,
\\0\ (\mo\ p^2)&\t{if}\ (\f{p}{267})=-1.\endcases
\endalign$$
Also,
$$\align &\f1{p^a}\sum_{k=0}^{p^a-1}\f{14151k+827}{(-300)^{3k}}\bi{2k}k^2\bi{3k}k
\\\eq& 827\l(\f {p^a}3\r)+\l(\f{p^{a-1}}3\r)\f{13}{150}p^2B_{p-2}\l(\f13\r)\ (\mo\ p^3)
\endalign$$
for any positive integer $a$.
\endproclaim
\Remark\ 2.10. It is known (cf. [CC]) that
$$\sum_{k=0}^\infty\f{615k+53}{(-48)^{3k}}\bi{2k}{k}^2\bi{3k}{k}=\f{96\sqrt3}{\pi}$$
and
$$\sum_{k=0}^\infty\f{14151k+827}{(-300)^{3k}}\bi{2k}{k}^2\bi{4k}{2k}=\f{1500\sqrt3}{\pi}.$$
\medskip

\heading{3. Using the polynomials $S_n(x)=\sum_{k=0}^n\bi nk^4x^k$}\endheading

In this section we give many conjectures on Problem 1.1 involving
the polynomials
$$S_n(x):=\sum_{k=0}^n\bi nk^4x^k\quad (n=0,1,2,\ldots).$$

\proclaim{Conjecture 3.1} Let $p$ be an odd prime. Then
$$\sum_{k=0}^{p-1}S_k(-2)\eq\cases 4x^2-2p\pmod{p^2}&\t{if}\ p=x^2+y^2\ (2\nmid x\ \&\ 2\mid y),
\\0\pmod{p^2}&\t{if}\ p\eq3\pmod 4.\endcases$$
And
$$\sum_{k=0}^{p-1}(3k+2)S_k(-2)\eq\f p2\l(1+3\l(\f{-1}p\r)\r)\pmod {p^2}.$$
Moreover,
$$\f1n\sum_{k=0}^{n-1}(6k+4)S_k(-2)\in\Z\qquad\t{for all}\ n=1,2,3,\ldots.$$
\endproclaim

\proclaim{Conjecture 3.2} Let $p$ be an odd prime. Then
$$\align&\sum_{k=0}^{p-1}S_k(-4)\eq\sum_{k=0}^{p-1}S_k(-64)
\\\eq&\cases 4x^2-2p\pmod{p^2}&\t{if}\ p\eq1,9\pmod{20}\ \&\ p=x^2+5y^2\ (x,y\in\Z),
\\2x^2-2p\pmod{p^2}&\t{if}\ p\eq 3,7\pmod{20}\ \&\ 2p=x^2+5y^2\ (x,y\in\Z),
\\0\pmod{p^2}&\t{if}\ (\f{-5}p)=-1,\ \t{i.e.,}\ p\eq11,13,17,19\pmod{20}.\endcases
\endalign$$
And
$$\sum_{k=0}^{p-1}(8k+7)S_k(-64)\eq p\l(\f p{15}\r)\l(3+4\l(\f{-1}p\r)\r)\pmod {p^2}.$$
Moreover,
$$\f1n\sum_{k=0}^{n-1}(8k+7)S_k(-64)\in\Z\qquad\t{for all}\ n=1,2,3,\ldots.$$
\endproclaim

\proclaim{Conjecture 3.3} Let $p$ be an odd prime. Then
$$\align&\sum_{k=0}^{p-1}S_k(4)
\\\eq&\cases 4x^2-2p\pmod{p^2}&\t{if}\ p\eq1,7\pmod{24}\ \&\ p=x^2+6y^2\ (x,y\in\Z),
\\8x^2-2p\pmod{p^2}&\t{if}\ p\eq 5,11\pmod{24}\ \&\ p=2x^2+3y^2\ (x,y\in\Z),
\\0\pmod{p^2}&\t{if}\ (\f{-6}p)=-1,\ \t{i.e.,}\ p\eq13,17,19,23\pmod{24}.\endcases
\endalign$$
And
$$\sum_{k=0}^{p-1}(24k+17)S_k(4)\eq p\l(5+12\l(\f{2}p\r)\r)\pmod {p^2}.$$
Moreover,
$$\f1n\sum_{k=0}^{n-1}(24k+17)S_k(4)\in\Z\qquad\t{for all}\ n=1,2,3,\ldots.$$
\endproclaim

\proclaim{Conjecture 3.4} Let $p$ be an odd prime. Then
$$\align&\sum_{k=0}^{p-1}S_k(16)
\\\eq&\cases 4x^2-2p\pmod{p^2}&\t{if}\ p\eq1,9,11,19\pmod{40}\ \&\ p=x^2+10y^2,
\\8x^2-2p\pmod{p^2}&\t{if}\ p\eq 7,13,23,37\pmod{40}\ \&\ p=2x^2+5y^2,
\\0\pmod{p^2}&\t{if}\ (\f{-10}p)=-1.\endcases
\endalign$$
When $(\f{-10}p)=1$ we have
$$\sum_{k=0}^{p-1}(160k+129)S_k(16)\eq 80p\l(\f p5\r)\pmod {p^2}.$$
\endproclaim

\proclaim{Conjecture 3.5} Let $p$ be an odd prime. Then
$$\align&\sum_{k=0}^{p-1}S_k(1)\eq\sum_{k=0}^{p-1}S_k(-9)
\\\eq&\cases 4x^2-2p\pmod{p^2}&\t{if}\ p\eq1,4\pmod{15}\ \&\ p=x^2+15y^2,
\\12x^2-2p\pmod{p^2}&\t{if}\ p\eq 2,8\pmod{15}\ \&\ p=3x^2+5y^2,
\\0\pmod{p^2}&\t{if}\ (\f p{15})=-1.\endcases
\endalign$$
And
$$\align \sum_{k=0}^{p-1}(3k+2)S_k(1)\eq& 2p\l(\f p5\r)\pmod {p^2},
\\\sum_{k=0}^{p-1}(5k+4)S_k(-9)\eq&\f p2\l(\f p3\r)\l(3+5\l(\f p{15}\r)\r)\pmod{p^2}.
\endalign$$
Moreover,
$$\f1{2n}\sum_{k=0}^{n-1}(3k+2)S_k(1)\in\Z
\ \  \t{and}\ \ \f1{2n}\sum_{k=0}^{n-1}(5k+4)S_k(-9)\in\Z$$
for all $n=1,2,3,\ldots$.
\endproclaim

\proclaim{Conjecture 3.6} Let $p$ be an odd prime. Then
$$\align&\sum_{k=0}^{p-1}S_k(36)
\\\eq&\cases 4x^2-2p\pmod{p^2}&\t{if}\ (\f 2p)=(\f p3)=(\f p5)=1\ \&\ p=x^2+30y^2,
\\12x^2-2p\pmod{p^2}&\t{if}\ (\f p3)=1,\ (\f 2p)=(\f p5)=-1\ \&\ p=3x^2+10y^2,
\\8x^2-2p\pmod{p^2}&\t{if}\ (\f 2p)=1,\ (\f p3)=(\f p5)=-1\ \&\ p=2x^2+15y^2,
\\2p-20x^2\pmod{p^2}&\t{if}\ (\f p5)=1,\ (\f 2p)=(\f p3)=-1\ \&\ p=5x^2+6y^2,
\\0\pmod{p^2}&\t{if}\ (\f{-30}p)=-1.\endcases
\endalign$$
And
$$\sum_{k=0}^{p-1}(8k+7)S_k(36)\eq p\l(\f p{15}\r)\l(3+4\l(\f{-6}p\r)\r)\pmod {p^2}.$$
We also have
$$\f1{n}\sum_{k=0}^{n-1}(8k+7)S_k(36)\in\Z
\qquad\t{for all}\ n=1,2,3,\ldots.$$
\endproclaim

\proclaim{Conjecture 3.7} Let $p$ be an odd prime. Then
$$\align&\sum_{k=0}^{p-1}S_k(196)
\\\eq&\cases 4x^2-2p\pmod{p^2}&\t{if}\ (\f 2p)=(\f p5)=(\f p7)=1\ \&\ p=x^2+70y^2,
\\8x^2-2p\pmod{p^2}&\t{if}\ (\f p7)=1,\ (\f 2p)=(\f p5)=-1\ \&\ p=2x^2+35y^2,
\\2p-20x^2\pmod{p^2}&\t{if}\ (\f p5)=1,\ (\f 2p)=(\f p7)=-1\ \&\ p=5x^2+14y^2,
\\28x^2-2p\pmod{p^2}&\t{if}\ (\f 2p)=1,\ (\f p5)=(\f p7)=-1\ \&\ p=7x^2+10y^2,
\\0\pmod{p^2}&\t{if}\ (\f{-70}p)=-1.\endcases
\endalign$$
And
$$\sum_{k=0}^{p-1}(120k+109)S_k(196)\eq p\l(\f p{7}\r)\l(49+60\l(\f{-14}p\r)\r)\pmod {p^2}.$$
We also have
$$\f1{n}\sum_{k=0}^{n-1}(120k+109)S_k(196)\in\Z
\qquad\t{for all}\ n=1,2,3,\ldots.$$
\endproclaim

\proclaim{Conjecture 3.8} Let $p$ be an odd prime. Then
$$\align&\sum_{k=0}^{p-1}S_k(-324)
\\\eq&\cases 4x^2-2p\pmod{p^2}&\t{if}\ (\f {-1}p)=(\f p5)=(\f p{17})=1\ \&\ p=x^2+85y^2,
\\2x^2-2p\pmod{p^2}&\t{if}\ (\f {p}{17})=1,\ (\f {-1}p)=(\f p5)=-1\ \&\ 2p=x^2+85y^2,
\\2p-20x^2\pmod{p^2}&\t{if}\ (\f {-1}p)=1,\ (\f p5)=(\f p{17})=-1\ \&\ p=5x^2+17y^2,
\\2p-10x^2\pmod{p^2}&\t{if}\ (\f p{5})=1,\ (\f {-1}p)=(\f p{17})=-1\ \&\ 2p=5x^2+17y^2,
\\0\pmod{p^2}&\t{if}\ (\f{-85}p)=-1.\endcases
\endalign$$
Provided $p>3$ we have
$$\sum_{k=0}^{p-1}(34k+31)S_k(-324)\eq p\l(\f p5\r)\l(17+14\l(\f{-1}p\r)\r)\pmod {p^2}.$$
Moreover,
$$\f1{n}\sum_{k=0}^{n-1}(34k+31)S_k(-324)\in\Z
\qquad\t{for all}\ n=1,2,3,\ldots.$$
\endproclaim

\proclaim{Conjecture 3.9} Let $p$ be an odd prime. Then
$$\align&\sum_{k=0}^{p-1}S_k(1296)
\\\eq&\cases 4x^2-2p\pmod{p^2}&\t{if}\ (\f {-2}p)=(\f p5)=(\f p{13})=1\ \&\ p=x^2+130y^2,
\\8x^2-2p\pmod{p^2}&\t{if}\ (\f {-2}p)=1,\ (\f p5)=(\f p{13})=-1\ \&\ p=2x^2+65y^2,
\\2p-20x^2\pmod{p^2}&\t{if}\ (\f p5)=1,\ (\f {-2}p)=(\f p{13})=-1\ \&\ p=5x^2+26y^2,
\\2p-40x^2\pmod{p^2}&\t{if}\ (\f p{13})=1,\ (\f {-2}p)=(\f p5)=-1\ \&\ p=10x^2+13y^2,
\\0\pmod{p^2}&\t{if}\ (\f{-130}p)=-1.\endcases
\endalign$$
Provided $p>3$ we have
$$\sum_{k=0}^{p-1}(130k+121)S_k(1296)\eq p\l(\f{-2} p\r)\l(56+65\l(\f{-26}p\r)\r)\pmod {p^2}.$$
Moreover,
$$\f1{n}\sum_{k=0}^{n-1}(130k+121)S_k(1296)\in\Z
\qquad\t{for all}\ n=1,2,3,\ldots.$$
\endproclaim

\proclaim{Conjecture 3.10} Let $p$ be an odd prime. Then
$$\align&\sum_{k=0}^{p-1}S_k(5776)
\\\eq&\cases 4x^2-2p\pmod{p^2}&\t{if}\ (\f {2}p)=(\f p5)=(\f p{19})=1\ \&\ p=x^2+190y^2,
\\8x^2-2p\pmod{p^2}&\t{if}\ (\f {2}p)=1,\ (\f p5)=(\f p{19})=-1\ \&\ p=2x^2+95y^2,
\\2p-20x^2\pmod{p^2}&\t{if}\ (\f p{19})=1,\ (\f {2}p)=(\f p{5})=-1\ \&\ p=5x^2+38y^2,
\\2p-40x^2\pmod{p^2}&\t{if}\ (\f p{5})=1,\ (\f {2}p)=(\f p{19})=-1\ \&\ p=10x^2+19y^2,
\\0\pmod{p^2}&\t{if}\ (\f{-190}p)=-1.\endcases
\endalign$$
And
$$\sum_{k=0}^{p-1}(816k+769)S_k(5776)\eq p\l(\f p{95}\r)\l(361+408\l(\f{p}{19}\r)\r)\pmod {p^2}.$$
Moreover,
$$\f1{n}\sum_{k=0}^{n-1}(816k+769)S_k(5776)\in\Z
\qquad\t{for all}\ n=1,2,3,\ldots.$$
\endproclaim

\proclaim{Conjecture 3.11} Let $p$ be an odd prime. Then
$$\align\sum_{k=0}^{p-1}S_k(12)
\eq\cases4x^2-2p\pmod{p^2}&\t{if}\ p\eq1\pmod{12}\ \&\ p=x^2+y^2\ (3\nmid x),
\\(\f{xy}3)4xy\pmod{p^2}&\t{if}\ p\eq5\pmod{12}\ \&\ p=x^2+y^2,
\\0\pmod{p^2}&\t{if}\ p\eq3\pmod 4.\endcases
\endalign$$
And
$$\sum_{k=0}^{p-1}(4k+3)S_k(12)\eq p\l(1+2\l(\f 3p\r)\r)\pmod {p^2}.$$
Moreover,
$$\f1{n}\sum_{k=0}^{n-1}(4k+3)S_k(12)\in\Z
\qquad\t{for all}\ n=1,2,3,\ldots.$$
\endproclaim

\proclaim{Conjecture 3.12} Let $p$ be an odd prime. Then
$$\align&\sum_{k=0}^{p-1}S_k(-20)
\\\eq&\cases4x^2-2p\pmod{p^2}&\t{if}\ p\eq1,9\pmod{20}\ \&\ p=x^2+y^2\ (5\nmid x),
\\4xy\pmod{p^2}&\t{if}\ p\eq13,17\pmod{20}\ \&\  p=x^2+y^2\ (5\mid x-y),
\\0\pmod{p^2}&\t{if}\ p\eq3\pmod 4.\endcases
\endalign$$
And
$$\sum_{k=0}^{p-1}(6k+5)S_k(-20)\eq p\l(\f{-1}p\r)\l(2+3\l(\f {-5}p\r)\r)\pmod {p^2}.$$
Moreover,
$$\f1{n}\sum_{k=0}^{n-1}(6k+5)S_k(-20)\in\Z
\qquad\t{for all}\ n=1,2,3,\ldots.$$
\endproclaim

In 2005 Yifan Yang [Y] found the interesting identity
$$\sum_{k=0}^\infty\f{4k+1}{36^k}\sum_{j=0}^k\bi kj^4=\f{18}{\pi\sqrt{15}}.$$
Motivated by this we give the following easy result.

\proclaim{Theorem 3.1} Let $a,b\in\Z$ and $m\in\Z\sm\{0\}$. Let $p$ be an odd prime not dividing $m$. Then
$$\sum_{k=0}^{p-1}\f{ak+b}{m^k}\sum_{j=0}^k\bi kj^4\eq-\sum_{k=0}^{p-1}(2ak+2a-b)S_k(m)\pmod p.$$
In particular,
$$\sum_{k=0}^{p-1}\f{\sum_{j=0}^k\bi kj^4}{m^k}\eq\sum_{k=0}^{p-1}S_k(m)\pmod{p}.$$
\endproclaim
\Proof. Let
$$S:=\sum_{n=0}^{p-1}(2an+2a-b)S_n(m).$$ Then
$$\align S=&\sum_{k=0}^{p-1}m^k\sum_{n=k}^{p-1}(2a(n+1)-b)\bi nk^4
\\=&\sum_{k=0}^{p-1}m^k\sum_{j=0}^{p-1-k}(2a(k+j+1)-b)\bi{k+j}j^4
\\=&\sum_{k=0}^{p-1}m^k\sum_{j=0}^{p-1-k}(2a(k+j+1)-b)\bi{-k-1}j^4
\\\eq&\sum_{k=0}^{p-1}m^k\sum_{j=0}^{p-1-k}(2a(k+j+1)-b)\bi{p-1-k}j^4\pmod p
\endalign$$
and hence
$$\align S\eq&\sum_{k=0}^{p-1}m^{p-1-k}\sum_{j=0}^{k}(2a(p-k+j)-b)\bi{k}j^4
\\\eq&-\sum_{k=0}^{p-1}\f1{m^k}\sum_{l=0}^k(2al+b)\bi kl^4\pmod p.
\endalign$$
Note that
$$2\sum_{l=0}^kl\bi kl^4=\sum_{j=0}^k\l(j\bi kj^4+(k-j)\bi k{k-j}^4\r)=k\sum_{j=0}^k\bi kj^4.$$
So we have
$$S\eq-\sum_{k=0}^{p-1}\f{ak+b}{m^k}\sum_{j=0}^k\bi kj^4\pmod p.$$
This completes the proof. \qed

\proclaim{Conjecture 3.13} Let $m$ be among
$$1,\,-2,\,\pm4,\,-9,\,12,\,16,\,-20,\,36,\,-64,\,196,\,-324,\,1296,\,5776.$$
For any odd prime $p$ not dividing $m$ we have
$$\sum_{k=0}^{p-1}\f{\sum_{j=0}^k\bi kj^4}{m^k}\eq\sum_{k=0}^{p-1}S_k(m)\pmod{p^2}.$$
\endproclaim

Now we give two more conjectures involving $\sum_{k=0}^n\bi nk^3x^k$.

\proclaim{Conjecture 3.14} Let $p$ be an odd prime. Then
$$\align&\sum_{n=0}^{p-1}\f{\bi{4n}{2n}}{16^n}\sum_{k=0}^n\bi nk^3
\\\eq&\cases(\f{-1}p)(4x^2-2p)\pmod{p^2}&\t{if}\ p\eq1,3\pmod8\ \&\ p=x^2+2y^2,
\\0\pmod{p^2}&\t{if}\ p\eq5,7\pmod8.\endcases
\endalign$$
When $p\eq5,7\pmod8$, we have
$$\sum_{k=0}^{p-1}(16n+5)\f{\bi{4n}{2n}}{16^n}\sum_{k=0}^n\bi nk^3\eq0\pmod{p^2}.$$
\endproclaim

\proclaim{Conjecture 3.15} {\rm (i)} Let $p>5$ be a prime with
$p\not=11$. Then
$$\align&\sum_{n=0}^{p-1}\f{\bi{2n}n}{(-121)^n}\sum_{k=0}^n\bi nk^3(-5)^{3k}
\\\eq&\cases 4x^2-2p\pmod{p^2}&\t{if}\ p\eq1,4\pmod{15}\ \&\ p=x^2+15y^2,
\\2p-12x^2\pmod{p^2}&\t{if}\ p\eq2,8\pmod{15}\ \&\ p=3x^2+15y^2,
\\0\pmod{p^2}&\t{if}\ p\eq7,11,13,14\pmod{15}.\endcases
\endalign$$

 {\rm (ii)} For any prime $p>5$ with $p\not=41$, we have
$$\align&\sum_{n=0}^{p-1}\f{\bi{2n}n}{4100^n}\sum_{k=0}^n\bi nk^32^{12k}
\\\eq&\cases x^2-2p\pmod{p^2}&\t{if}\ (\f p3)=(\f p{41})=1\ \&\ 4p=x^2+123y^2,
\\2p-3x^2\pmod{p^2}&\t{if}\ (\f p3)=(\f p{41})=-1\ \&\ 4p=3x^2+41y^2,
\\0\pmod{p^2}&\t{if}\ (\f p{123})=-1.\endcases
\endalign$$

{\rm (iii)}  Let $p>3$ be a prime with $p\not=53,89$. Then
$$\align&\sum_{n=0}^{p-1}\f{\bi{2n}n}{1000004^n}\sum_{k=0}^n\bi nk^310^{6k}
\\\eq&\cases x^2-2p\pmod{p^2}&\t{if}\ (\f p3)=(\f p{89})=1\ \&\ 4p=x^2+267y^2,
\\2p-3x^2\pmod{p^2}&\t{if}\ (\f p3)=(\f p{89})=-1\ \&\ 4p=3x^2+89y^2,
\\0\pmod{p^2}&\t{if}\ (\f p{267})=-1.\endcases
\endalign$$
\endproclaim

\heading{4. Using the function $F_n(x)=\sum_{k=0}^n\bi nk^3\bi{2k}kx^{-k}$}\endheading

For $n=0,1,2,\ldots$ we define
$$F_n(x):=\sum_{k=0}^n\bi nk^3\bi{2k}kx^{-k}.$$

\proclaim{Conjecture 4.1} Let $p$ be an odd prime. Then
$$\sum_{k=0}^{p-1}(-1)^kF_k(-2)
\eq\cases4x^2-2p\pmod{p^2}&\t{if}\ p\eq1\pmod{12}\ \&\ p=x^2+y^2\ (3\nmid x),
\\(\f{xy}3)4xy\pmod{p^2}&\t{if}\ p\eq5\pmod{12}\ \&\ p=x^2+y^2,
\\0\pmod{p^2}&\t{if}\ p\eq3\pmod4.\endcases$$
And
$$\sum_{k=0}^{p-1}(3k+2)(-1)^kF_k(-2)\eq\f p2\l(\f{-1}p\r)\l(3\l(\f p3\r)+1\r)\pmod{p^2}.$$
\endproclaim

\proclaim{Conjecture 4.2} Let $p$ be an odd prime. Then
$$\sum_{k=0}^{p-1}(-1)^kF_k(2)
\eq\cases4x^2-2p\pmod{p^2}&\t{if}\ p\eq1,9\pmod{20}\ \&\ p=x^2+5y^2,
\\2x^2-2p\pmod{p^2}&\t{if}\ p\eq3,7\pmod{20}\ \&\ 2p=x^2+y^2,
\\0\pmod{p^2}&\t{if}\  p\eq11,13,17,19\pmod{20}.\endcases$$
And
$$\sum_{k=0}^{p-1}(10k+7)(-1)^kF_k(2)\eq p\l(5\l(\f {-1}p\r)+2\r)\pmod{p^2}.$$
\endproclaim

\proclaim{Conjecture 4.3} Let $p>3$ be a prime. Then
$$\sum_{k=0}^{p-1}(-1)^kF_k(4)
\eq\cases4x^2-2p\pmod{p^2}&\t{if}\ p\eq1,7\pmod{24}\ \&\ p=x^2+6y^2,
\\8x^2-2p\pmod{p^2}&\t{if}\ p\eq5,11\pmod{24}\ \&\ p=2x^2+3y^2,
\\0\pmod{p^2}&\t{if}\ p\eq13,17,19,23\pmod{24}.\endcases$$
When $p\eq1,5,7,11\pmod{24}$, we have
$$\sum_{k=0}^{p-1}(72k+47)(-1)^kF_k(4)\eq 0\pmod{p^2}.$$
\endproclaim

\proclaim{Conjecture 4.4} Let $p>3$ be a prime. Then
$$\align&\sum_{k=0}^{p-1}(-1)^kF_k(12)
\\\eq&\cases(\f p3)(4x^2-2p)\pmod{p^2}&\t{if}\ p\eq1,3\pmod{8}\ \&\ p=x^2+2y^2,
\\0\pmod{p^2}&\t{if}\  p\eq5,7\pmod{8}.\endcases\endalign$$
And
$$\sum_{k=0}^{p-1}(5k+3)(-1)^kF_k(12)\eq\f p2\l(\f{-2}p\r)\l(5\l(\f p3\r)+1\r)\pmod{p^2}.$$
\endproclaim

\proclaim{Conjecture 4.5} Let $p\not=2,7$ be a prime. Then
$$\align&\sum_{k=0}^{p-1}(-1)^kF_k(-14)
\\\eq&\cases4x^2-2p\pmod{p^2}&\t{if}\ (\f{-1}p)=(\f p3)=(\f p7)=1\ \&\ p=x^2+21y^2,
\\2x^2-2p\pmod{p^2}&\t{if}\ (\f p7)=1,\,(\f{-1}p)=(\f p3)=-1\ \&\ 2p=x^2+21y^2,
\\2p-12x^2\pmod{p^2}&\t{if}\ (\f{p}3)=1,\,(\f {-1}p)=(\f p7)=-1\ \&\ p=3x^2+7y^2,
\\2p-6x^2\pmod{p^2}&\t{if}\ (\f{-1}p)=1,\,(\f p3)=(\f p7)=-1\ \&\ 2p=3x^2+7y^2,
\\0\pmod{p^2}&\t{if}\ (\f{-21}p)=-1.\endcases
\endalign$$
And
$$\sum_{k=0}^{p-1}(18k+11)(-1)^kF_k(-14)\eq p\l(\f p7\r)\l(2+9\l(\f{-1}p\r)\r)\pmod{p^2}.$$
\endproclaim

\proclaim{Conjecture 4.6} Let $p>5$ be a prime. Then
$$\align&\sum_{k=0}^{p-1}(-1)^kF_k(40)
\\\eq&\cases4x^2-2p\pmod{p^2}&\t{if}\ (\f{2}p)=(\f p3)=(\f p5)=1\ \&\ p=x^2+30y^2,
\\8x^2-2p\pmod{p^2}&\t{if}\ (\f{2}p)=1,\,(\f p3)=(\f p5)=-1\ \&\ p=2x^2+15y^2,
\\2p-12x^2\pmod{p^2}&\t{if}\ (\f{p}3)=1,\,(\f 2p)=(\f p5)=-1\ \&\ p=3x^2+10y^2,
\\20x^2-2p\pmod{p^2}&\t{if}\ (\f{p}5)=1,\,(\f 2p)=(\f p3)=-1\ \&\ p=5x^2+6y^2,
\\0\pmod{p^2}&\t{if}\ (\f{-30}p)=-1.\endcases
\endalign$$
And
$$\sum_{k=0}^{p-1}(12k+7)(-1)^kF_k(40)\eq p\l(\f{10}p\r)\l(6+\l(\f{-6}p\r)\r)\pmod{p^2}.$$
\endproclaim

\proclaim{Conjecture 4.7} Let $p$ be an odd prime. We have
$$\align&\sum_{k=0}^{p-1}(-1)^kF_k(-50)
\\\eq&\cases4x^2-2p\pmod{p^2}&\t{if}\ (\f{-1}p)=(\f p3)=(\f p{11})=1\ \&\ p=x^2+33y^2,
\\2x^2-2p\pmod{p^2}&\t{if}\ (\f{-1}p)=1,\,(\f p3)=(\f p{11})=-1\ \&\ 2p=x^2+33y^2,
\\2p-12x^2\pmod{p^2}&\t{if}\ (\f{p}{11})=1,\,(\f {-1}p)=(\f p3)=-1\ \&\ p=3x^2+11y^2,
\\2p-6x^2\pmod{p^2}&\t{if}\ (\f{p}3)=1,\,(\f {-1}p)=(\f p{11})=-1\ \&\ 2p=3x^2+11y^2,
\\0\pmod{p^2}&\t{if}\ (\f{-33}p)=-1.\endcases
\endalign$$
If $p\not=5$, then
$$\sum_{k=0}^{p-1}(99k+58)(-1)^kF_k(-50)\eq \f p2\l(\f{-1}p\r)\l(99\l(\f{p}3\r)+17\r)\pmod{p^2}.$$
\endproclaim

\proclaim{Conjecture 4.8} Let $p$ be an odd prime. We have
$$\align&\sum_{k=0}^{p-1}(-1)^kF_k(112)
\\\eq&\cases4x^2-2p\pmod{p^2}&\t{if}\ (\f{-2}p)=(\f p3)=(\f p{7})=1\ \&\ p=x^2+42y^2,
\\8x^2-2p\pmod{p^2}&\t{if}\ (\f{p}7)=1,\,(\f {-2}p)=(\f p{3})=-1\ \&\ p=2x^2+21y^2,
\\2p-12x^2\pmod{p^2}&\t{if}\ (\f{-2}{p})=1,\,(\f p3)=(\f p7)=-1\ \&\ p=3x^2+14y^2,
\\2p-24x^2\pmod{p^2}&\t{if}\ (\f{p}3)=1,\,(\f {-2}p)=(\f p{7})=-1\ \&\ p=6x^2+7y^2,
\\0\pmod{p^2}&\t{if}\ (\f{-42}p)=-1.\endcases
\endalign$$
If $p\not=7$, then
$$\sum_{k=0}^{p-1}(180k+103)(-1)^kF_k(112)\eq p\l(\f{p}7\r)\l(90\l(\f{p}3\r)+13\r)\pmod{p^2}.$$
\endproclaim

\proclaim{Conjecture 4.9} Let $p$ be an odd prime. We have
$$\align&\sum_{k=0}^{p-1}(-1)^kF_k(-338)
\\\eq&\cases4x^2-2p\pmod{p^2}&\t{if}\ (\f{-1}p)=(\f p3)=(\f p{19})=1\ \&\ p=x^2+57y^2,
\\2x^2-2p\pmod{p^2}&\t{if}\ (\f{-1}p)=1,\,(\f p3)=(\f p{19})=-1\ \&\ 2p=x^2+57y^2,
\\2p-12x^2\pmod{p^2}&\t{if}\ (\f{p}3)=1,\,(\f {-1}p)=(\f p{19})=-1\ \&\ p=3x^2+19y^2,
\\2p-6x^2\pmod{p^2}&\t{if}\ (\f{p}{19})=1,\,(\f {-1}p)=(\f p{3})=-1\ \&\ 2p=3x^2+19y^2,
\\0\pmod{p^2}&\t{if}\ (\f{-57}p)=-1.\endcases
\endalign$$
If $p\not=13$, then
$$\sum_{k=0}^{p-1}(855k+482)(-1)^kF_k(-338)\eq \f p2\l(\f{-1}p\r)\l(855\l(\f{p}{19}\r)+109\r)\pmod{p^2}.$$
\endproclaim

\proclaim{Conjecture 4.10} Let $p>5$ be a prime with $p\not=13$. Then
$$\align&\sum_{k=0}^{p-1}(-1)^kF_k(1300)
\\\eq&\cases 4x^2-2p\pmod{p^2}&\t{if}\ (\f 2p)=(\f p3)=(\f p{13})=1\ \&\ p=x^2+78y^2,
\\8x^2-2p\pmod{p^2}&\t{if}\ (\f 2p)=1,\ (\f p3)=(\f p{13})=-1\ \&\ p=2x^2+39y^2,
\\2p-12x^2\pmod{p^2}&\t{if}\ (\f p{13})=1,\ (\f 2p)=(\f p3)=-1\ \&\ p=3x^2+26y^2,
\\2p-24x^2\pmod{p^2}&\t{if}\ (\f p3)=1,\ (\f 2p)=(\f p{13})=-1\ \&\ p=6x^2+13y^2,
\\0\pmod{p^2}&\t{if}\ (\f{-78}p)=-1.\endcases
\endalign$$
And
$$\sum_{k=0}^{p-1}(204k+113)(-1)^kF_k(1300)\eq p\l(\f p{39}\r)\l(102\l(\f p3\r)+11\r)\pmod{p^2}.$$
\endproclaim

\proclaim{Conjecture 4.11} Let $p>3$ be a prime with $p\not=31$. We have
$$\align&\sum_{k=0}^{p-1}(-1)^kF_k(-3038)
\\\eq&\cases 4x^2-2p\pmod{p^2}&\t{if}\ (\f {-1}p)=(\f p3)=(\f p{31})=1\ \&\ p=x^2+93y^2,
\\8x^2-2p\pmod{p^2}&\t{if}\ (\f p{31})=1,\ (\f {-1}p)=(\f p{3})=-1\ \&\ 2p=x^2+93y^2,
\\2p-12x^2\pmod{p^2}&\t{if}\ (\f p{3})=1,\ (\f {-1}p)=(\f p{31})=-1\ \&\ p=3x^2+31y^2,
\\2p-6x^2\pmod{p^2}&\t{if}\ (\f {-1}p)=1,\ (\f p3)=(\f p{31})=-1\ \&\ 2p=3x^2+31y^2,
\\0\pmod{p^2}&\t{if}\ (\f{-93}p)=-1.\endcases
\endalign$$
If $p\not=7$, then
$$\sum_{k=0}^{p-1}(1170k+643)(-1)^kF_k(-3038)\eq p\l(\f p{31}\r)\l(585\l(\f {-1}p\r)+58\r)\pmod{p^2}.$$
\endproclaim

\proclaim{Conjecture 4.12} Let $p>3$ be a prime.  We have
$$\align&\sum_{k=0}^{p-1}(-1)^kF_k(4900)
\\\eq&\cases 4x^2-2p\pmod{p^2}&\t{if}\ (\f {2}p)=(\f p3)=(\f p{17})=1\ \&\ p=x^2+102y^2,
\\8x^2-2p\pmod{p^2}&\t{if}\ (\f p{17})=1,\ (\f {2}p)=(\f p{3})=-1\ \&\ p=2x^2+51y^2,
\\2p-12x^2\pmod{p^2}&\t{if}\ (\f p{3})=1,\ (\f {2}p)=(\f p{17})=-1\ \&\ p=3x^2+34y^2,
\\2p-24x^2\pmod{p^2}&\t{if}\ (\f {2}p)=1,\ (\f p3)=(\f p{17})=-1\ \&\ p=6x^2+17y^2,
\\0\pmod{p^2}&\t{if}\ (\f{-102}p)=-1.\endcases
\endalign$$
If $p\not=7$, then
$$\sum_{k=0}^{p-1}(561k+307)(-1)^kF_k(4900)\eq \f p2\l(\f {-6}p\r)\l(561\l(\f {p}{51}\r)+53\r)\pmod{p^2}.$$
\endproclaim

\heading{5. Using the function $G_n(x)=\sum_{k=0}^n\bi nk^2\bi{2k}k\bi{2n-2k}{n-k}x^{-k}$}
\endheading

For $n\in\N$ we define
$$G_n(x):=\sum_{k=0}^n\bi nk^2\bi{2k}k\bi{2n-2k}{n-k}x^{-k}.$$
Those integers $G_n(1)\ (n=0,1,2,\ldots)$ are called Domb numbers.
In 2004 H. H. Chan, S. H. Chan and Z. Liu [CCL] proved that
$$\sum_{k=0}^\infty\f{5k+1}{64^k}G_k(1)=\f 8{\sqrt3\,\pi}.$$
Motivated by his work with Mahler measures nad new transformation
formulas for ${}_5F_4$ series, M. D. Rogers [Ro] discovered that
$$\sum_{k=0}^\infty\f{3k+1}{32^k}(-1)^kG_k(1)=\f 2{\pi},$$
which was independently showed by H. H. Chan and H. Verrill [CV].

\proclaim{Conjecture 5.1} Let $p>3$ be a prime.  We have
$$\align&\sum_{k=0}^{p-1}G_k(1)\eq\sum_{k=0}^{p-1}G_k(-5)\eq\sum_{k=0}^{p-1}\f{G_k(1)}{64^k}
\\\eq&\cases 4x^2-2p\pmod{p^2}&\t{if}\ p\eq1,4\pmod{15}\ \&\ p=x^2+15y^2,
\\2p-12x^2\pmod{p^2}&\t{if}\ p\eq2,8\pmod{p^2}\ \&\ p=3x^2+5y^2,
\\0\pmod{p^2}&\t{if}\ (\f{-15}p)=-1.\endcases
\endalign$$
If $p\not=5$, then
$$\sum_{k=0}^{p-1}(5k+4)G_k(1)\eq 4p\l(\f p3\r)\pmod{p^3}$$
and $$\sum_{k=0}^{p-1}(21k+16)G_k(-5)\eq 16p\l(\f p5\r)\pmod{p^2}.$$
Also,
$$\f1{4n}\sum_{k=0}^{n-1}(5k+4)G_k(1)\in\Z\qquad\t{for all}\ n=1,2,3,\ldots.$$
\endproclaim
\proclaim{Conjecture 5.2} Let $p>3$ be a prime.  We have
$$\align&\sum_{k=0}^{p-1}\f{G_k(1)}{(-2)^k}\eq\sum_{k=0}^{p-1}\f{G_k(1)}{4^k}\eq\sum_{k=0}^{p-1}\f{G_k(1)}{16^k}
\eq\sum_{k=0}^{p-1}\f{G_k(1)}{(-32)^k}\eq\l(\f{-1}p\r)\sum_{k=0}^{p-1}G_k(6)
\\\eq&\cases 4x^2-2p\pmod{p^2}&\t{if}\ p\eq1\pmod{3}\ \&\ p=x^2+3y^2\ (x,y\in\Z),
\\0\pmod{p^2}&\t{if}\ p\eq2\pmod{3}.\endcases
\endalign$$
Also,
$$\align\sum_{k=0}^{p-1}(3k+2)\f{G_k(1)}{(-2)^k}\eq& 2p\l(\f {-1}p\r)\pmod{p^3},
\\\sum_{k=0}^{p-1}(3k+1)\f{G_k(1)}{(-32)^k}\eq& p\l(\f {-1}p\r)+p^3E_{p-3}\pmod{p^4},
\\\sum_{k=0}^{p-1}(15k+11)G_k(6)\eq& p\l(10\l(\f3p\r)+1\r)\pmod{p^2}.
\endalign$$
If $p\eq1\pmod3$, then
$$\sum_{k=0}^{p-1}(3k+2)\f{G_k(1)}{4^k}\eq\sum_{k=0}^{p-1}(3k+1)\f{G_k(1)}{16^k}\eq0\pmod{p^2}.$$
\endproclaim

\proclaim{Conjecture 5.3} Let $p>3$ be a prime.  Then
$$\sum_{k=0}^{p-1}\f{G_k(1)}{8^k}
\eq\cases 4x^2-2p\pmod{p^2}&\t{if}\ p\eq1,3\pmod{8}\ \&\ p=x^2+2y^2,
\\0\pmod{p^2}&\t{if}\ p\eq5,7\pmod{8}.\endcases$$
And
$$\sum_{k=0}^{p-1}\f{G_k(1)}{(-8)^k}
\eq\cases4x^2-2p\pmod{p^2}&\t{if}\ p\eq1,7\pmod{24}\ \&\ p=x^2+6y^2,
\\8x^2-2p\pmod{p^2}&\t{if}\ p\eq5,11\pmod{24}\ \&\ p=2x^2+3y^2,
\\0\pmod{p^2}&\t{if}\ (\f{-6}p)=-1.
\endcases$$
Also,
$$\sum_{k=0}^{p-1}(2k+1)\f{G_k(1)}{8^k}\eq p\ (\mo\ p^4)\
\ \t{and}\ \ \sum_{k=0}^{p-1}(2k+1)\f{G_k(1)}{(-8)^k}\eq p\l(\f p3\r)\ (\mo\ p^3).$$
\endproclaim

\proclaim{Conjecture 5.4} Let $p$ be an odd prime. We have
$$\align \sum_{k=0}^{p-1}G_k(4)\eq&\sum_{k=0}^{p-1}(-1)^kF_k(4)\pmod{p^2},
\\\sum_{k=0}^{p-1}G_k(-12)\eq&\l(\f 3p\r)\sum_{k=0}^{p-1}(-1)^kF_k(-14)\pmod{p^2}\ \ (p\not=3,7),
\\\sum_{k=0}^{p-1}G_k(36)\eq&\l(\f 2p\r)\sum_{k=0}^{p-1}(-1)^kF_k(40)\pmod{p^2}\ \ (p>5),
\\\sum_{k=0}^{P-1}G_k(-44)\eq&\sum_{k=0}^{p-1}(-1)^kF_k(-50)\pmod{p^2}\ \ (p\not=5,11),
\\\sum_{k=0}^{p-1}G_k(100)\eq&\sum_{k=0}^{p-1}(-1)^kF_k(112)\pmod{p^2}\ \ (p\not=5,7),
\\\sum_{k=0}^{p-1}G_k(-300)\eq&\l(\f{-1}p\r)\sum_{k=0}^{p-1}(-1)^kF_k(-338)\pmod{p^2}\ \ (p\not=3,5,13),
\\\sum_{k=0}^{p-1}G_k(1156)\eq&\sum_{k=0}^{p-1}(-1)^kF_k(1300)\pmod{p^2}\ \ (p\not=5,13,17),
\\\sum_{k=0}^{p-1}G_k(-2700)\eq&\l(\f{3}p\r)\sum_{k=0}^{p-1}(-1)^kF_k(-3038)\pmod{p^2}\ \ (p>7,\, p\not=31),
\\\sum_{k=0}^{p-1}G_k(4356)\eq&\l(\f{-6}p\r)\sum_{k=0}^{p-1}(-1)^kF_k(4900)\pmod{p^2}\ \ (p>11).
\endalign$$
\endproclaim

\heading{6. Using $a_n(x)=\sum_{k=0}^n\bi nk^2\bi{n+k}kx^k$}\endheading

For $n\in\N$ define
$$a_n(x):=\sum_{k=0}^n\bi nk^2\bi{n+k}kx^k.$$
Those numbers $a_n(1)\ (n=0,1,2,\ldots)$ first appeared in Ap\'ery's
proof of the irrationality of $\zeta(2)$ (see [Ap] and [P]). We
observe the new identity
$$\bi{2n}na_n(1)=\sum_{k=0}^n\bi
nk^2\bi{2k}n\bi{n+2k}n.\tag6.1$$ If we set $u_n=\sum_{k=0}^n\bi
nk^2\bi{2k}n\bi{n+2k}n/\bi{2n}n$ or $u_n=a_n(1)$ for $n\in\N$, then
$u_0=1$ and $u_1=3$, and by applying the Zeilberger algorithm (cf.
[PWZ]) via {\tt Mathematica} (version 7) we get the recurrence
relation
$$(n+2)^2u_{n+2}=(11n^2+33n+25)u_{n+1}+(n+1)^2u_n\ \ (n=0,1,2,\ldots) .$$
Thus (6.1) holds by induction.

We find $\sum_{k=0}^{p-1}\bi{2k}ka_k(1)/{m^k}$ mod $p^2$ related to
the representation $4p=x^2+dy^2$ with
$$\align(m,d)=&(-3,15),(4,11),(18,1),(-28,35),(36,19),(72,10),
\\&(147,15),(-828,115),(-15228,235).
\endalign$$

\proclaim{Conjecture 6.1}  Let $p$ be an odd prime. Then
$$\sum_{k=0}^{p-1}\f{\bi{2k}ka_k(1)}{4^k}
\eq\cases x^2-2p\pmod{p^2}&\t{if}\ (\f p{11})=1\ \&\ 4p=x^2+11y^2,
\\0\pmod{p^2}&\t{if}\ (\f{p}{11})=-1.\endcases$$
And
$$\sum_{k=0}^{p-1}\f{22k+9}{4^k}\bi{2k}ka_k(1)\eq 9p\pmod{p^2}.$$
\endproclaim

\proclaim{Conjecture 6.2} Let $p>3$ be a prime. Then
$$\align\sum_{k=0}^{p-1}\f{\bi{2k}ka_k(1)}{36^k}
\eq\cases x^2-2p\pmod{p^2}&\t{if}\ (\f p{19})=1\ \&\ 4p=x^2+19y^2,
\\0\pmod{p^2}&\t{if}\ (\f{p}{19})=-1.\endcases
\endalign$$
And
$$\sum_{k=0}^{p-1}\f{38k+13}{36^k}\bi{2k}ka_k(1)\eq 13p\pmod{p^2}.$$
\endproclaim

\proclaim{Conjecture 6.3}  Let $p\not=2,7$ be a prime. Then
$$\sum_{k=0}^{p-1}\f{\bi{2k}ka_k(1)}{(-28)^k}
\eq\cases x^2-2p\pmod{p^2}&\t{if}\ (\f p5)=(\f p{7})=1\ \&\ 4p=x^2+35y^2,
\\2p-5x^2\pmod{p^2}&\t{if}\ (\f p5)=(\f p{7})=-1\ \&\ 4p=5x^2+7y^2,
\\0\pmod{p^2}&\t{if}\ (\f p{35})=-1.\endcases$$
And
$$\sum_{k=0}^{p-1}\f{10k+3}{(-28)^k}\bi{2k}ka_k(1)\eq 3p\l(\f p{7}\r)\pmod{p^2}.$$
\endproclaim

\proclaim{Conjecture 6.4} Let $p\not=2,3,23$ be a prime. Then
$$\sum_{k=0}^{p-1}\f{\bi{2k}ka_k(1)}{(-828)^k}
\eq\cases x^2-2p\pmod{p^2}&\t{if}\ (\f p5)=(\f p{23})=1\ \&\ 4p=x^2+115y^2,
\\2p-5x^2\pmod{p^2}&\t{if}\ (\f p5)=(\f p{23})=-1\ \&\ 4p=5x^2+23y^2,
\\0\pmod{p^2}&\t{if}\ (\f{p}{115})=-1.\endcases$$
And
$$\sum_{k=0}^{p-1}\f{190k+29}{(-828)^k}\bi{2k}ka_k(1)\eq 29p\l(\f p{23}\r)\pmod{p^2}.$$
\endproclaim

\proclaim{Conjecture 6.5} Let $p\not=5,47$ be an odd prime. Then
$$\sum_{k=0}^{p-1}\f{\bi{2k}ka_k(1)}{(-15228)^k}
\eq\cases x^2-2p\pmod{p^2}&\t{if}\ (\f p5)=(\f p{47})=1\ \&\ 4p=x^2+235y^2,
\\2p-5x^2\pmod{p^2}&\t{if}\ (\f p5)=(\f p{47})=-1\ \&\ 4p=5x^2+47y^2,
\\0\pmod{p^2}&\t{if}\ (\f{p}{235})=-1.\endcases$$
When $p\not=3$ we have
$$\sum_{k=0}^{p-1}\f{682k+71}{(-15228)^k}\bi{2k}ka_k(1)\eq 71p\l(\f p{47}\r)\pmod{p^2}.$$
\endproclaim

\proclaim{Conjecture 6.6} Let $p>5$ be a prime. Then
$$\align&\sum_{k=0}^{p-1}\f{\bi{2k}ka_k(-3)}{4^k}\eq \sum_{k=0}^{p-1}\f{\bi{2k}ka_k(9)}{100^k}
\\\eq&\cases x^2-2p\pmod{p^2}&\t{if}\ p\eq1\pmod{3}\ \&\ 4p=x^2+27y^2,
\\0\pmod{p^2}&\t{if}\  p\eq2\pmod{3}.
\endcases
\endalign$$
Also,
$$\sum_{k=0}^{p-1}(30k+13)\f{\bi{2k}ka_k(-3)}{4^k}\eq p\l(3+10\l(\f{p}3\r)\r)\pmod{p^2}$$
and
$$\sum_{k=0}^{p-1}(66k+23)\f{\bi{2k}ka_k(9)}{100^k}\eq 23p\pmod{p^2}.$$
\endproclaim

\proclaim{Conjecture 6.7} Let $p>3$ be a prime. Then
$$\align\sum_{k=0}^{p-1}\f{\bi{2k}ka_k(-27)}{26^{2k}}
\eq\cases (\f p3)(x^2-2p)\pmod{p^2}&\t{if}\ (\f p{11})=1\ \&\ 4p=x^2+11y^2,
\\0\pmod{p^2}&\t{if}\  (\f p{11})=-1.
\endcases
\endalign$$
When $p\not=7$, we have
$$\align\sum_{k=0}^{p-1}\f{\bi{2k}ka_k(27)}{28^{2k}}
\eq\cases (\f p3)(4x^2-2p)\pmod{p^2}&\t{if}\ p\eq1,3\pmod8\ \&\ p=x^2+2y^2,
\\0\pmod{p^2}&\t{if}\  p\eq5,7\pmod{8}.
\endcases
\endalign$$
Also,
$$\sum_{k=0}^{\infty}(114k+31)\f{\bi{2k}ka_k(-27)}{26^{2k}}=\f{338\sqrt3}{11\pi}\ \t{and}\
\sum_{k=0}^{\infty}(930k+143)\f{\bi{2k}ka_k(27)}{28^{2k}}=\f{980\sqrt3}{\pi}.$$
\endproclaim

For $n\in\N$ we define
$$a_n^*(x):=\sum_{k=0}^n\bi nk^2\bi{n+k}kx^{n-k}=x^na_n\l(\f1x\r).$$

\proclaim{Conjecture 6.8} Let $p>3$ be a prime. Then
$$\align&\l(\f{-1}p\r)\sum_{k=0}^{p-1}\f{\bi{2k}ka_k^*(-8)}{96^k}
\\\eq&\cases 4x^2-2p&\t{if}\ p\eq1,7\pmod{24}\ \&\ p=x^2+6y^2,
\\8x^2-2p\pmod{p^2}&\t{if}\ p\eq5,11\pmod{24}\ \&\ p=2x^2+3y^2,
\\0\pmod{p^2}&\t{if}\ (\f{-6}p)=-1.\endcases
\endalign$$
And
$$\sum_{k=0}^{p-1}(13k+4)\f{\bi{2k}ka_k^*(-8)}{96^k}\eq 4p\l(\f{-2}p\r)\pmod{p^2}.$$
We also have
$$\sum_{k=0}^\infty(13k+4)\f{\bi{2k}ka_k^*(-8)}{96^k}=\f{9\sqrt2}{2\pi}.$$
\endproclaim

\proclaim{Conjecture 6.9} Let $p>3$ be a prime. Then
$$\sum_{k=0}^{p-1}\f{\bi{2k}ka_k^*(-32)}{1152^k}
\eq\cases (\f 2p)(4x^2-2p)&\t{if}\ p\eq1\pmod{3}\ \&\ p=x^2+3y^2,
\\0\pmod{p^2}&\t{if}\ p\eq2\pmod3.\endcases$$
And
$$\sum_{k=0}^{p-1}(290k+61)\f{\bi{2k}ka_k^*(-32)}{1152^k}\eq p\l(\f2p\r)\l(6+55\l(\f{-1}p\r)\r)\pmod{p^2}.$$
We also have
$$\sum_{k=0}^\infty(290k+61)\f{\bi{2k}ka_k^*(-32)}{1152^k}=\f{99\sqrt2}{\pi}.$$
\endproclaim
\proclaim{Conjecture 6.10} Let $p>5$ be a prime. Then
$$\align&\l(\f{-1}p\r)\sum_{k=0}^{p-1}\f{\bi{2k}ka_k^*(64)}{3840^k}
\\\eq&\cases 4x^2-2p\pmod{p^2}&\t{if}\ p\eq1,4\pmod{15}\ \&\ p=x^2+15y^2,
\\2p-12x^2\pmod{p^2}&\t{if}\ p\eq2,8\pmod{15}\ \&\ p=3x^2+5y^2,
\\0\pmod{p^2}&\t{if}\ (\f{p}{15})=-1.\endcases
\endalign$$
And
$$\sum_{k=0}^{p-1}(962k+137)\f{\bi{2k}ka_k^*(64)}{3840^k}\eq p\l(\f{-5}p\r)\l(147-10\l(\f p3\r)\r)\pmod{p^2}.$$
We also have
$$\sum_{k=0}^\infty(962k+137)\f{\bi{2k}ka_k^*(64)}{3840^k}=\f{252\sqrt5}{\pi}.$$
\endproclaim

\heading{7. Miscellaneous things}\endheading

 In this section we give some miscellaneous conjectures related to Problem 1.2.

\proclaim{Conjecture 7.1} Let $p>3$ be a prime. Then
$$\align&\sum_{n=0}^{p-1}\f1{9^n}\sum_{k=0}^n\bi {-1/3}k^2\bi{-2/3}{n-k}^2
\\\eq&\cases x^2-2p\pmod{p^2}&\t{if}\ p\eq1\pmod{3}\ \&\ 4p=x^2+27y^2,
\\0\pmod{p^2}&\t{if}\ p\eq2\pmod3.\endcases
\endalign$$
And
$$\sum_{n=0}^{p-1}\f{8n+1}{9^n}\sum_{k=0}^n\bi {-1/3}k^2\bi{-2/3}{n-k}^2
\eq p\l(\f{p}3\r)\pmod{p^3}.$$
We also have
$$\sum_{n=0}^\infty\f{8n+1}{9^n}\sum_{k=0}^n\bi {-1/3}k^2\bi{-2/3}{n-k}^2=\f{3\sqrt3}{\pi}.$$
\endproclaim

\proclaim{Conjecture 7.2} Let $p>3$ be a prime. Then
$$\align&\sum_{n=0}^{p-1}\f1{2^n}\sum_{k=0}^n\bi {-1/3}k\bi{-2/3}{n-k}\bi{-1/6}k\bi{-5/6}{n-k}
\\\eq&\cases 4x^2-2p\pmod{p^2}&\t{if}\ p\eq1,7\pmod{24}\ \&\ p=x^2+6y^2,
\\2p-8x^2\pmod{p^2}&\t{if}\ p\eq5,11\pmod{24}\ \&\ p=2x^2+3y^2,
\\0\pmod{p^2}&\t{if}\ (\f{-6}{p})=-1.\endcases
\endalign$$
And
$$\sum_{n=0}^{p-1}\f{3n-1}{2^n}\sum_{k=0}^n\bi {-1/3}k\bi{-2/3}{n-k}\bi{-1/6}k\bi{-5/6}{n-k}
\eq -p\l(\f{-6}p\r)\pmod{p^2}.$$
We also have
$$\sum_{n=0}^\infty\f{3n-1}{2^n}\sum_{k=0}^n\bi {-1/3}k\bi{-2/3}{n-k}\bi{-1/6}k\bi{-5/6}{n-k}=\f{3\sqrt6}{2\pi}.$$
\endproclaim

\proclaim{Conjecture 7.3} Let $p$ be an odd prime. Then
$$\align&\l(\f{-1}p\r)\sum_{n=0}^{p-1}\f1{64^n}\sum_{k=0}^n\bi nk\bi{2k}n\bi{2k}k\bi{2n-2k}{n-k}3^{2k-n}
\\\eq&\cases 4x^2-2p\pmod{p^2}&\t{if}\ p\eq1,4\pmod{15}\ \&\ p=x^2+15y^2,
\\2p-12x^2\pmod{p^2}&\t{if}\ p\eq2,8\pmod{15}\ \&\ p=3x^2+5y^2,
\\0\pmod{p^2}&\t{if}\ (\f{p}{15})=-1.\endcases
\endalign$$
And
$$\align&\sum_{n=0}^{p-1}\f{21n+1}{64^n}\sum_{k=0}^n\bi nk\bi{2k}n\bi{2k}k\bi{2n-2k}{n-k}3^{2k-n}
\\&\qquad\eq p\l(\f{-1}p\r)\l(4-3\l(\f p3\r)\r)\pmod{p^2}.
\endalign$$
We also have
$$\sum_{n=0}^\infty\f{21n+1}{64^n}\sum_{k=0}^n\bi nk\bi{2k}n\bi{2k}k\bi{2n-2k}{n-k}3^{2k-n}=\f{64}{\pi}.$$
\endproclaim

\proclaim{Conjecture 7.4} Let $p$ be an odd prime. Then
$$\align&\sum_{n=0}^{p-1}\f1{(-64)^n}\sum_{k=0}^n\bi nk^2\bi{n+k}k\bi{2k}n
\\\eq&
\cases x^2-2p\pmod{p^2}&\t{if}\ (\f p7)=(\f p{13})=1\ \&\
4p=x^2+91y^2,
\\2p-7x^2\pmod{p^2}&\t{if}\ (\f p7)=(\f p{13})=-1\ \&\ 4p=7x^2+13y^2,
\\0\pmod{p^2}&\t{if}\ (\f p{91})=-1.\endcases
\endalign$$
And
$$\sum_{n=0}^{p-1}(39n+10)\sum_{k=0}^n\bi nk^2\bi{n+k}k\bi{2k}n
\eq10p\l(\f p7\r)\pmod{p^2}.$$
\endproclaim

\proclaim{Conjecture 7.5} Let $p>3$ be a prime. Then
$$\align&\sum_{n=0}^{p-1}\f1{216^n}\sum_{k=0}^n\bi nk^2\bi{n+k}k\bi{2k}n14^{2k-n}
\\\eq&
\cases x^2-2p\pmod{p^2}&\t{if}\ (\f p{67})=1\ \&\
4p=x^2+67y^2,
\\0\pmod{p^2}&\t{if}\ (\f p{67})=-1.\endcases
\endalign$$
And
$$\sum_{n=0}^{p-1}(33165n+11546)\sum_{k=0}^n\bi nk^2\bi{n+k}k\bi{2k}n14^{2k-n}
\eq11546p\pmod{p^2}.$$
\endproclaim

\proclaim{Conjecture 7.6} Let $p$ be an odd prime. Then
$$\align&\sum_{n=0}^{p-1}\f1{(-8)^n}\sum_{k=0}^n\bi nk^2\bi{n+k}k\bi{2k}n26^{2k-n}
\\\eq&\cases x^2-2p\pmod{p^2}&\t{if}\ (\f p5)=(\f p{23})=1\ \&\ 4p=x^2+115y^2,
\\5x^2-2p\pmod{p^2}&\t{if}\ (\f p5)=(\f p{23})=-1\ \&\ 4p=5x^2+23y^2,
\\0\pmod{p^2}&\t{if}\ (\f p{115})=-1.\endcases
\endalign$$
And
$$\align&\sum_{n=0}^{p-1}\f{6555n+3062}{(-8)^n}\sum_{k=0}^n\bi nk^2\bi{n+k}k\bi{2k}n26^{2k-n}
\\&\qquad\eq 2p\l(220+1311\l(\f p{23}\r)\r)\pmod{p^2}.
\endalign$$
\endproclaim

\proclaim{Conjecture 7.7} Let $p$ be an odd prime. Then
$$\align&\sum_{n=0}^{p-1}\f1{(-8)^n}\sum_{k=0}^n\bi nk^2\bi{3k}n\bi{3(n-k)}n
\\\eq&\cases x^2-2p\pmod{p^2}&\t{if}\ (\f p5)=(\f p{7})=1\ \&\ 4p=x^2+35y^2,
\\2p-5x^2\pmod{p^2}&\t{if}\ (\f p5)=(\f p{7})=-1\ \&\ 4p=5x^2+7y^2,
\\0\pmod{p^2}&\t{if}\ (\f p{35})=-1.\endcases
\endalign$$
And
$$\sum_{n=0}^{p-1}\f{35n+18}{(-8)^n}\sum_{k=0}^n\bi nk^2\bi{3k}n\bi{3(n-k)}n
\eq 9p\l(7-5\l(\f p{5}\r)\r)\pmod{p^2}.$$
\endproclaim

\proclaim{Conjecture 7.8} Let $p>3$ be a prime. When $p\not=11,17$, we have
$$\align&\sum_{n=0}^{p-1}\f{\bi{2n}n}{198^{2n}}\sum_{k=0}^n\bi nk^2\bi{2k}k
\\\eq&\cases 4x^2-2p\pmod{p^2}&\t{if}\ (\f {2}p)=(\f p3)=(\f p{17})=1\ \&\ p=x^2+102y^2,
\\2p-8x^2\pmod{p^2}&\t{if}\ (\f p{17})=1,\ (\f {2}p)=(\f p{3})=-1\ \&\ p=2x^2+51y^2,
\\12x^2-2p\pmod{p^2}&\t{if}\ (\f p{3})=1,\ (\f {2}p)=(\f p{17})=-1\ \&\ p=3x^2+34y^2,
\\2p-24x^2\pmod{p^2}&\t{if}\ (\f {2}p)=1,\ (\f p3)=(\f p{17})=-1\ \&\ p=6x^2+17y^2,
\\0\pmod{p^2}&\t{if}\ (\f{-102}p)=-1.\endcases
\endalign$$
When $p\not=23,59$, we have
$$\align&\sum_{n=0}^{p-1}\f{\bi{2n}n}{(-1123596)^n}\sum_{k=0}^n\bi nk^2\bi{2k}k
\\\eq&\cases 4x^2-2p\pmod{p^2}&\t{if}\ (\f {-1}p)=(\f p{3})=(\f p{59})=1\ \&\ p=x^2+177y^2,
\\2p-2x^2\pmod{p^2}&\t{if}\ (\f{-1}p)=1,\ (\f p3)=(\f p{59})=-1\ \&\ 2p=x^2+177y^2,
\\12x^2-2p\pmod{p^2}&\t{if}\ (\f{p}{59})=1,\ (\f {-1}p)=(\f p3)=-1\ \&\ p=3x^2+59y^2,
\\2p-6x^2\pmod{p^2}&\t{if}\ (\f{p}{3})=1,\ (\f {-1}p)=(\f p{59})=-1\ \&\ 2p=3x^2+59y^2,
\\0\pmod{p^2}&\t{if}\ (\f {-177}p)=-1.\endcases
\endalign$$

\endproclaim

\proclaim{Conjecture 7.9} Let $p>5$ be a prime. Then
$$\align&\sum_{n=0}^{p-1}\f{\bi{2n}n}{324^n}\sum_{k=0}^n\bi nk^2\bi{2k}k(-20)^k
\\\eq&\cases 4x^2-2p\pmod{p^2}&\t{if}\ p\eq1,9\ (\mo\ 20)\ \&\ p=x^2+y^2\ (5\nmid x),
\\4xy\pmod{p^2}&\t{if}\ p\eq13,17\ (\mo\ 20)\ \&\ p=x^2+y^2\ (5\mid x-y),
\\0\pmod{p^2}&\t{if}\ p\eq3\ (\mo\ 4).\endcases
\endalign$$
And
$$\sum_{n=0}^{p-1}\f{16n+5}{324^n}\bi{2n}n\sum_{k=0}^n\bi nk^2\bi{2k}k(-20)^k\eq\f p{25}\l(112\l(\f {-1}p\r)+13\r)\pmod{p^2}.$$
We also have
$$\sum_{n=0}^\infty\f{16n+5}{324^n}\bi{2n}n\sum_{k=0}^n\bi nk^2\bi{2k}k(-20)^k=\f{189}{25\pi}.$$
\endproclaim

\proclaim{Conjecture 7.10} Let $p>5$ be a prime. Then
$$\align&\sum_{n=0}^{p-1}\f{\bi{2n}n}{400^n}\sum_{k=0}^n\bi nk\bi{n+2k}{2k}\bi{2k}k196^{n-k}
\\\eq&\cases x^2-2p\pmod{p^2}&\t{if}\ (\f 2p)=(\f p{11})=1\ \&\ p=x^2+22y^2,
\\8x^2-2p\pmod{p^2}&\t{if}\ (\f 2p)=(\f p{11})=-1\ \&\ p=2x^2+11y^2,
\\0\pmod{p^2}&\t{if}\ (\f {-22}p)=-1.\endcases
\endalign$$
And
$$\sum_{n=0}^{p-1}\f{33n+19}{400^n}\bi{2n}n\sum_{k=0}^n\bi nk\bi{n+2k}{2k}\bi{2k}k196^{n-k}\eq19p\pmod{p^2}.$$
\endproclaim

For $n\in\N$ we define
$$f_n^-(x):=\sum_{k=0}^n\bi nk^2\bi{2k}n(-1)^kx^{2k-n}.$$

\proclaim{Conjecture 7.11} Let $p>5$ be a prime. Then
$$\align&\l(\f p3\r)\sum_{k=0}^{p-1}\f{\bi{2k}kf_k^{-}(8)}{480^k}
\\\eq&\cases4x^2-2p\pmod{p^2}&\t{if}\ p\eq1,9\pmod{20}\ \&\ p=x^2+5y^2,
\\2x^2-2p\pmod{p^2}&\t{if}\ p\eq3,7\pmod{20}\ \&\ 2p=x^2+5y^2,
\\0\pmod{p^2}&\t{if}\ p\eq 11,13,17,19\pmod{20},
\endcases
\endalign$$
and
$$\align&\sum_{k=0}^{p-1}\f{\bi{2k}kf_k^-(18)}{5760^k}
\\\eq&\cases 4x^2-2p\pmod{p^2}&\t{if}\ (\f {-1}p)=(\f p5)=(\f p{17})=1\ \&\ p=x^2+85y^2,
\\2p-2x^2\pmod{p^2}&\t{if}\ (\f {p}{17})=1,\ (\f {-1}p)=(\f p5)=-1\ \&\ 2p=x^2+85y^2,
\\20x^2-2p\pmod{p^2}&\t{if}\ (\f {-1}p)=1,\ (\f p5)=(\f p{17})=-1\ \&\ p=5x^2+17y^2,
\\2p-10x^2\pmod{p^2}&\t{if}\ (\f p{5})=1,\ (\f {-1}p)=(\f p{17})=-1\ \&\ 2p=5x^2+17y^2,
\\0\pmod{p^2}&\t{if}\ (\f{-85}p)=-1.\endcases
\endalign$$
Also,
$$\sum_{k=0}^{p-1}\f{1054k+233}{480^k}\bi{2k}kf_k^{-}(8)
\eq p\l(\f{-1}p\r)\l(221+12\l(\f p{15}\r)\r)\pmod{p^2}$$
and
$$\sum_{k=0}^{p-1}\f{224434k+32849}{5760^k}\bi{2k}kf_k^-(18)\eq p\l(32305\l(\f{-1}p\r)+544\l(\f p5\r)\r)\pmod{p^2}.$$
We also have
$$\sum_{k=0}^\infty\f{1054k+233}{480^k}\bi{2k}kf_k^-(8)=\f{520}{\pi}$$
and
$$\sum_{k=0}^\infty\f{224434k+32849}{5760^k}\bi{2k}kf_k^-(18)=\f{93600}{\pi}.$$
\endproclaim

For $n=0,1,2,\ldots$, we define
$$f_n^+(x):=\sum_{k=0}^n\bi nk^2\bi{2k}nx^{2k-n}.$$
Recall that $f_n^+(1)=\sum_{k=0}^n\bi nk^3$ by an identity of V. Strehl (see also D. Zagier [Z]).

\proclaim{Conjecture 7.12} {\rm (i)} Let $p>3$ be a prime. Then
$$\align\sum_{k=0}^{p-1}\f{\bi{2k}kf_k^+(9)}{144^k}
\eq\cases4x^2-2p\pmod{p^2}&\t{if}\ (\f{-1}p)=(\f p{13})=1\ \&\ p=x^2+13y^2,
\\2x^2-2p\pmod{p^2}&\t{if}\  (\f{-1}p)=(\f p{13})=-1\ \&\ 2p=x^2+13y^2,
\\0\pmod{p^2}&\t{if}\ (\f{-13}p)=-1.\endcases
\endalign$$
And $$\sum_{k=0}^{p-1}\f{65k+22}{144^k}\bi{2k}kf_k^+(9)
\eq 22p\pmod{p^2}.$$

{\rm (ii)} Let $p\not=5$ be an odd prime. Then
$$\align&\sum_{k=0}^{p-1}\f{\bi{2k}kT_k^2(13,13)}{1300^k}
\\\eq&\cases 4x^2-2p\ (\mo\ p^2)&\t{if}\ (\f{-1}p)=(\f{13}p)=1\ \&\ p=x^2+13y^2,
\\2x^2-2p\ (\mo\ p^2)&\t{if}\ (\f{-1}p)=(\f {13}p)=-1\ \&\ 2p=x^2+13y^2,
\\0\ (\mo\ p^2)&\t{if}\ (\f{-13}p)=-1.\endcases\endalign$$
And
$$\sum_{k=0}^{p-1}\f{312k+91}{1300^k}\bi{2k}kT_k^2(13,13)\eq 91p\pmod{p^2}.$$
\endproclaim

\proclaim{Conjecture 7.13} Let $p>3$ be a prime. Then
$$\align&\sum_{k=0}^{p-1}\f{\bi{2k}kf^+_k(-7)}{16^k}\eq\sum_{k=0}^{p-1}\f{\bi{2k}kT_k^2(3,-3)}{(-108)^k}
\eq\l(\f p3\r)\sum_{n=0}^{p-1}\f{\bi{2n}n}{(-108)^n}\sum_{k=0}^n\bi nk^2\bi{2k}k
\\\eq&\cases 4x^2-2p\pmod{p^2}&\t{if}\ (\f{-1}p)=(\f{p}3)=(\f p{7})=1\ \&\ p=x^2+21y^2,
\\12x^2-2p\pmod{p^2}&\t{if}\ (\f{-1}p)=(\f p7)=-1,\ (\f p3)=1\ \&\ p=3x^2+7y^2,
\\2x^2-2p\pmod{p^2}&\t{if}\ (\f {-1}p)=(\f p{3})=-1,\ (\f p{7})=1\ \&\ 2p=x^2+21y^2,
\\6x^2-2p\pmod{p^2}&\t{if}\ (\f {-1}p)=1,\ (\f p{3})=(\f p7)=-1,\ \&\ 2p=3x^2+7y^2,
\\0\pmod{p^2}&\t{if}\ (\f{-21}p)=-1.\endcases
\endalign$$
We also have
$$\sum_{k=0}^{p-1}\f{63k+26}{16^k}\bi{2k}kf_k^+(-7)\eq p\l(5+21\l(\f p7\r)\r)\pmod{p^2}$$
and
$$\sum_{k=0}^{p-1}\f{56k+19}{(-108)^k}\bi{2k}kT_{k}^2(3,-3)
\eq \f p2\l(21\l(\f{p}7\r)+17\r)\pmod{p^2}.$$
\endproclaim

\proclaim{Conjecture 7.14} Let $p\not=43$ be an odd prime. Then
$$\align&\sum_{k=0}^{p-1}\f{\bi{2k}kf_k^+(175)}{(-29584)^k}
\\\eq&\cases4x^2-2p\pmod{p^2}&\t{if}\ (\f{-1}p)=(\f p7)=(\f p{19})=1\ \&\ p=x^2+133y^2,
\\2x^2-2p\pmod{p^2}&\t{if}\ (\f p{7})=1,\ (\f{-1}p)=(\f p{19})=-1\ \&\ 2p=x^2+133y^2,
\\2p-28x^2\pmod{p^2}&\t{if}\ (\f{p}{19})=1,\ (\f{-1}p)=(\f p{7})=-1\ \&\ p=7x^2+19y^2,
\\2p-14x^2\pmod{p^2}&\t{if}\ (\f {-1}p)=1,\ (\f{p}7)=(\f p{19})=-1\ \&\ 2p=7x^2+19y^2,
\\0\pmod{p^2}&\t{if}\ (\f{-133}p)=-1.\endcases
\endalign$$
And $$\sum_{k=0}^{p-1}\f{8851815k+1356374}{(-29584)^k}\bi{2k}kf_k^+(175)
\eq p\l(1300495\l(\f p7\r)+55879\r)\pmod{p^2}.$$
Moreover,
$$\sum_{k=0}^\infty\f{8851815k+1356374}{(-29584)^k}\bi{2k}kf_k^+(175)=\f{1349770\sqrt7}{\pi}.$$
\endproclaim

We also find $\sum_{k=0}^{p-1}\bi{2k}kf^+_k(x)/{m^k}$ mod $p^2$ (with $p>3$ a prime not dividing $m\in\Z\sm\{0\}$)
related to the representation $4p=x^2+dy^2$
with
$$\align (m,x,d)=&(2,16,6),(4,-1,7),(4,80,10),(5,64,10),(6,16,2),(6,240,30),
\\&(7,-16,21),(14,289,7),(14,2800,70),(15,96,30),(20,16,30),
\\&(21,576,42),(24,400,2),(25,-16,33),(25,384,2),
\\&(36,336,42),(45,2304,70),(49,800,22),(50,784,22),
\\&(56,16,42),(144,720,70),(169,-16,57),(441,7056,37).
\endalign$$

Actually we have many other conjectural congruences and series for $1/\pi$, which cannot be listed here due to the limitation of the length of
this survey.

\medskip

\Ack. The author thanks the two referees for their helpful comments.

\enddocument